\documentclass[10pt]{article}
\usepackage{geometry}
\usepackage{amssymb} 
\usepackage{amsmath} 
\usepackage{stmaryrd} 
\usepackage{comment} 
\usepackage{color}
\usepackage{graphicx}

\usepackage{environ} 

\NewEnviron{CommentsFEM}
{ 
}             

\newcommand{\V}[1]{\underline{#1}}
\newcommand{\M}[1]{\underline{\underline{#1}}}

\newcommand{\Cola}{  \color{black}  }
\newcommand{\Colb}{  \color{black}  }

\newcommand{\Def}{{ \Cola \, := \, }}

\newcommand{\Domain}{{ \Cola \Omega  }}
\newcommand{\Dimension}{{ \Cola d  }}
\newcommand{\DDomain}{{ \Cola \partial \Omega }}
\newcommand{\DSurface}{{ \Cola d\Gamma }}
\newcommand{\DVolume}{{ \Cola d \Omega }}

\newcommand{\DDomainDirichlet}{{ \Cola \partial \Omega^\textrm{w} }}
\newcommand{\DDomainNeumann}{{ \Cola \partial \Omega^\textrm{t} }}

\newcommand{\Point}{{ \Cola M }}
\newcommand{\Pos}{{ \Cola \V{x} }}
\newcommand{\ReferenceFrame}{{ \Cola \mathcal{R} }}
\newcommand{\GlobalOrigin}{{ \Cola \V{O}_\mathcal{R} }}
\newcommand{\UnitVector}[1]{{ \Cola \V{e}_{#1} }}

\newcommand{\GlobalCoordinate}[1]{{ \Cola x_{#1} }}

\newcommand{\NeumannBC}{{ \Cola \V{t} }}
\newcommand{\DirichletBC}{{ \Cola \V{w} }}
\newcommand{\VolumeForce}{{ \Cola \V{b} }}
\newcommand{\NeumannBCParameter}{{ \Cola \NeumannBC(\Parameter) }}
\newcommand{\DirichletBCParameter}{{ \Cola \DirichletBC(\Parameter) }}
\newcommand{\VolumeForceParameter}{{ \Cola \VolumeForce(\Parameter) }}

\newcommand{\Disp}{{ \Cola  \V{u} }}
\newcommand{\DispAlt}{{ \Cola  \V{v} }}
\newcommand{\DispFEMAlt}{{ \Cola  \DispAlt^\text{h} }}
\newcommand{\DispFEM}{{ \Cola \Disp^\text{h} }}

\newcommand{\DispLift}{{ \Cola \Disp^\text{p} }}

\newcommand{\DispFEMHomo}{{ \Cola \Disp^\text{h,0} }}
\newcommand{\DispFEMLift}{{ \Cola \Disp^\text{h,p} }}
\newcommand{\ErrorFEM}{{ \Cola \V{e}^\text{h} }}
\newcommand{\ErrorReduced}{{ \Cola \V{e}^\text{r} }}
\newcommand{\ErrorReducedApprox}{{ \Cola \V{\widetilde{e}}^\text{r} }}
\newcommand{\Stress}{{ \Cola \M{\sigma} }}
\newcommand{\StressTest}{{ \Cola \Stress^\starß }}
\newcommand{\StressAlt}{{ \Cola \Stress^\star }}
\newcommand{\StressHat}{{ \Colb \widehat{\Stress} }}
\newcommand{\StressFEM}{{ \Colb \Stress^\text{h} }}

\newcommand{\StressReducedHat}{{ \Colb \widehat{\Stress} }}
\newcommand{\StressFEMHomo}{{ \Colb \Stress^\text{h,0} }}
\newcommand{\StressFEMLift}{{ \Colb \Stress^\text{h,p} }}
\newcommand{\StressFEMReducedHomo}{{ \Colb \Stress^\text{r,0} }}

\newcommand{\StressReduced}{{ \Colb \Stress^\text{r} }}
\newcommand{\DispTest}{{ \Cola \Disp^\star }}

\newcommand{\StrainBasic}{{ \Cola \M{\epsilon} }}
\newcommand{\StrainBasicEngin}{{ \Cola \V{\epsilon} }}
\newcommand{\Strain}{{ \Colb \StrainBasic(\Disp) }}

\newcommand{\StrainTest}{{ \Colb \StrainBasic(\DispTest) }}

\newcommand{\StrainReduced}{{ \Colb \StrainBasic(\DispReduced) }}

\newcommand{\DispParameter}{{ \Colb \Disp(\Parameter) }}

\newcommand{\Hooke}{{ \Cola \underset{\widetilde{}}{D} }}
\newcommand{\HookeEngin}{{ \Cola \M{D} }}
\newcommand{\Compliance}{{ \Cola \underset{\widetilde{}}{C} }}
\newcommand{\Young}{{ \Cola E }}
\newcommand{\Poisson}{{ \Cola \nu }}
\newcommand{\LameFirst}{{ \Cola \lambda }}
\newcommand{\LameSecond}{{ \Cola G }}

\newcommand{\Sobolev}[1]{{  \Cola \mathcal{H}^{#1}  }}
\newcommand{\SpaceDisp}{{ \Cola \mathcal{U}(\Domain) }}
\newcommand{\SpaceDispCoarse}{{ \Cola \mathcal{\bar{U}}(\Domain) }}
\newcommand{\SpaceDispFluctuation}{{ \Cola \mathcal{\widetilde{U}}(\Domain) }}

\newcommand{\SpaceDispAd}{{ \Cola \mathcal{U}^\text{Ad}(\Domain) }}
\newcommand{\SpaceDispAdParameter}{{ \Cola \mathcal{U}^\text{Ad}(\Domain;\Parameter) }}
\newcommand{\SpaceDispFEM}{{ \Cola \mathcal{U}^\text{h}(\Domain) }}
\newcommand{\SpaceDispRigid}{{ \Cola \mathcal{U}^\text{rigid}(\Domain ) }}
\newcommand{\SpaceDispFEMMinusRigid}{{ \Cola \mathcal{\widetilde{U}}^\text{h}(\Domain) }}
\newcommand{\SpaceDispFEMAd}{{ \Cola \mathcal{U}^\text{h,Ad}(\Domain) }}
\newcommand{\SpaceDispFEMAdParameter}{{ \Cola \mathcal{U}^\text{h,Ad}(\Domain;\Parameter) }}

\newcommand{\SpaceDispFEMZero}{{ \Cola \mathcal{U}^\text{h,0}(\Domain) }}
\newcommand{\SpaceDispZero}{{ \Cola \mathcal{U}^{\text{Ad},0}(\Domain) }}
\newcommand{\SpaceDispRedZero}{{ \Cola \mathcal{U}^\text{r,0}(\Domain) }}
\newcommand{\SpaceDispReducedEnhancedZero}{{ \Cola \mathcal{U}^\text{2r,0}(\Domain) }}
\newcommand{\SpaceStress}{{ \Cola \mathcal{S}(\Domain) }}
\newcommand{\SpaceStressAd}{{ \Cola \mathcal{S}^{\text{Ad}}(\Domain) }}
\newcommand{\SpaceStressAdParameter}{{ \Cola \mathcal{S}^{\text{Ad}}(\Domain;\Parameter) }}

\newcommand{\SpaceStressFEMAd}{{ \Cola \mathcal{S}^{\text{h,Ad}}(\Domain) }}
\newcommand{\SpaceStressFEMAdParameter}{{ \Cola \mathcal{S}^{\text{h,Ad}}(\Domain;\Parameter) }}
\newcommand{\SpaceStressFEMZero}{{ \Cola \mathcal{S}^{\text{h,0}}(\Domain) }}
\newcommand{\SpaceStressReducedParameter}{{ \Cola \mathcal{S}^\text{r}(\Domain;\Parameter) }}
\newcommand{\SpaceStressReducedZero}{{ \Cola \mathcal{S}^\text{r,0}(\Domain) }}

\newcommand{\ParameterSpace}{{ \Cola \mathcal{P} }}
\newcommand{\ParameterSpaceDiscr}{{ \Cola \mathcal{\widetilde{P}} }}

\newcommand{\ShapeFEMI}[1]{{ \Cola N_{#1} }}
\newcommand{\NumberShape}{{ \Cola n_\text{n} }}
\newcommand{\NumberUnknowns}{{ \Cola n_\text{u} }}
\newcommand{\MatrixShapeFEM}{{ \Cola \M{N} }}
\newcommand{\MatrixDerivativeShapeFEM}{{ \Cola \M{B} }}
\newcommand{\StiffnessFEM}{{ \Cola \M{\bold{K}} }}
\newcommand{\ForceFEM}{{ \Cola \V{\bold{F}} }}
\newcommand{\MassFEM}{{ \Cola \M{\bold M} }}

\newcommand{\NumberBilinearAffine}{{ \Cola n_\text{a} }}
\newcommand{\NumberLinearAffine}{{ \Cola n_\text{l} }}
\newcommand{\NumberLinearLiftAffine}{{ \Cola n_\text{p} }}
\newcommand{\NumberNeumannBCAffine}{{ \Cola n_\text{t} }}
\newcommand{\NeumannBCAffineI}[1]{{ \Cola \bar{\NeumannBC}_{#1} }}
\newcommand{\CoeffNeumannBCAffineI}[1]{{ \Cola \gamma^\text{t}_{#1} }}
\newcommand{\NumberVolumeForceAffine}{{ \Cola n_\text{b} }}
\newcommand{\VolumeForceAffineI}[1]{{ \Cola \bar{\VolumeForce}_{#1} }}
\newcommand{\CoeffVolumeForceAffineI}[1]{{ \Cola \gamma^\text{b}_{#1} }}
\newcommand{\NumberHookeAffine}{{ \Cola n_\text{d} }}
\newcommand{\HookeAffineI}[1]{{ \Cola \underset{\widetilde{}}{\bar{D}_{#1}} }}
\newcommand{\CoeffHookeAffineI}[1]{{ \Cola \gamma^\text{d}_{#1} }}
\newcommand{\NumberComplianceAffine}{{ \Cola n_\text{c} }}
\newcommand{\ComplianceAffineI}[1]{{ \Cola \underset{\widetilde{}}{\bar{C}_{#1}} }}
\newcommand{\CoeffComplianceAffineI}[1]{{ \Cola \gamma^\text{c}_{#1} }}
\newcommand{\NumberDirichletBCAffine}{{ \Cola n_\text{w} }}
\newcommand{\DirichletBCAffineI}[1]{{ \Cola \bar{\DirichletBC}_{#1} }}
\newcommand{\CoeffDirichletBCAffineI}[1]{{ \Cola \gamma^\text{w}_{#1} }}

\newcommand{\FunctionAdStaticAffineI}[1]{{ \Cola \widetilde{\psi}_{#1} }}
\newcommand{\StressAdStaticAffineI}[1]{{ \Cola \bar{\Stress}_{#1}^\text{h,p} }}
\newcommand{\FunctionAdKinematicAffineI}[1]{{ \Cola \V{\psi}_{#1} }}
\newcommand{\Riesz}{{ \Cola \widetilde{\Disp}^\text{h} }}

\newcommand{\RieszLift}{{ \Cola \widetilde{\Disp}^\text{h,p} }}

\newcommand{\NumberSnapshot}{{ \Cola n_\text{s} }}
\newcommand{\Snapshot}{{ \Cola \Xi }}
\newcommand{\SnapshotZero}{{ \Cola \Snapshot^\text{0} }}
\newcommand{\SnapshotStress}{{ \Cola \widetilde{\Xi} }}
\newcommand{\SnapshotStressZero}{{ \Cola \SnapshotStress^\text{0} }}

\newcommand{\ParameterSnapshot}[1]{{ \Cola \Parameter^\text{s}_{#1} }}

\newcommand{\ReducedBasisI}[1]{{ \Cola \V{\phi}_{#1} }}
\newcommand{\ReducedBasisStressI}[1]{{ \Cola \M{\widetilde{\phi}}_{#1} }}
\newcommand{\ReducedCoeffI}[1]{{ \Cola \alpha_{#1} }}
\newcommand{\ReducedCoeffStressI}[1]{{ \Cola \widetilde{\alpha}_{#1} }}
\newcommand{\ReducedCoeff}{{ \Cola \V{\alpha} }}
\newcommand{\DispReduced}{{ \Cola \Disp^\text{r} }}
\newcommand{\DispReducedEnhanced}{{ \Cola \Disp^\text{2r} }}
\newcommand{\DispReducedParameter}{\Disp^\text{r}(\Parameter)}

\newcommand{\DispReducedHomo}{{ \Cola \Disp^\text{r,0} }}
\newcommand{\DispReducedEnhancedHomo}{{ \Cola \Disp^\text{2r,0} }}
\newcommand{\NumberReducedBasis}{{ \Cola n_\phi }}
\newcommand{\NumberReducedBasisEnriched}{{ \Cola n^\text{2r}_\phi }}
\newcommand{\NumberReducedBasisStress}{{ \Cola \tilde{n}_\phi }}

\newcommand{\CorrelationOperator }{{ \Cola \M{\bold H}^\text{s} }}
\newcommand{\CorrelationOperatorStress }{{ \Cola \M{\widetilde{\bold H}}^\text{s} }}

\newcommand{\IdentityTensor}{{ \Cola \M{I}_d }}

\newcommand{\Parameter}{{ \Cola \V{\mu} }}
\newcommand{\ParameterI}[1]{{ \Cola \mu_{#1} }}
\newcommand{\NumberParameter}{{ \Cola n_\mu }}
\newcommand{\ParameterZero}{{ \Cola \Parameter_\text{0} }}
\newcommand{\Output}{{ \Cola { \V{Q} } }}

\newcommand{\CRE}{{ \Cola  \nu^\text{up} }}
\newcommand{\LowerBound}{{ \Cola  \nu^\text{low} }}

\title{Certification of projection-based reduced order modelling in computational homogenisation by the Constitutive Relation Error.  \\
}

\author{P. Kerfriden$^{1\footnote{pierre.kerfriden(at)gmail.com}}$, J. J. R\'odenas$^{2}$ and S. P.-A. Bordas$^{1}$
\\ \\
$^{1}$ Cardiff University, School of Engineering, \\ The Parade, CF243AA Cardiff, United Kingdom \\
$^{2}$ Centro de Investigac\'ion de Tecnolog\'ia de Veh\'iculos (CITV) \\
Universitat Poli\`ecnica de Val\`encia, E-46022-Valencia, Spain \\ 
}
\date{}
  
\begin{document}

\maketitle

\begin{abstract}
In this paper, we propose upper and lower error bounding techniques for reduced order modelling applied to the computational homogenisation of random composites. The upper bound relies on the construction of a reduced model for the stress field. Upon ensuring that the reduced stress satisfies the equilibrium in the finite element sense, the desired bounding  property is obtained. The lower bound is obtained by defining a hierarchical enriched reduced model for the displacement. We show that the sharpness of both error estimates can be seamlessly controlled by adapting the parameters of the corresponding reduced order model.
\\

\noindent \textbf{keywords:} Model Order Reduction, Error Estimation, Computational Homogenisation, Proper Orthogonal Decomposition, Constitutive Relation Error
\end{abstract}


\section{Introduction}

Reduced order modelling is becoming an increasingly popular tool to solve parametrised or time-dependent problems (\textit{e.g.} \cite{legresleyalonso2001,meyermatthies2003,kunishvolkwein2003,rozzahuynh2008,ansallemfarhat2008,Nouy2010,kerfridenpassieux2011}). Such problems appear in a number of applications in solid mechanics, including the treatment of uncertainties, structural optimisation and multiscale modelling. Here, we are particularly interested in the case of nested computational homogenisation schemes for random heterogeneous materials (\textit{e.g.} \cite{zohdiwriggers2005,fish2006,geerskouznetsova2010}) with parametrised micro-structural properties. One of the numerical bottlenecks of such approaches is the numerical costs associated to the representative volume element (RVE). Indeed, these costs become prohibitive when the material hetorogeneities are parametrised in order to proceed to their identification or optimisation with respect to some macroscopic overall properties.


Reduced order modelling proposes to deliver a surrogate model for the solution to a  parametrised problem, whose evaluation should be inexpensive an accurate.  Such techniques consist of two phases. A training (or ``offline" ) stage and an evaluation (or ``online" ) stage. During the training stage, the parameter domain is explored, which provides training data that are used to build the surrogate model over the parameter domain. In the evaluation stage, the surrogate model is evaluated at a particular point of interest in the parameter domain. Reduced order modelling techniques differ in the way they explore the parameter domain, define and evaluate the surrogate model. One of the classical family of reduced order models is based on the response surface method. The solution is evaluated at certain points of the parameter domain, and interpolated using explicitly defined shape functions, for instance polynomials. A more advanced class of interpolation techniques such as Kriging, Moving Least-Squares or the Reduced Basis Method \cite{prudhommerovas2002,senveroy2006,rozzahuynh2008,hoangkhoo2013}, do not explicitly define the shape functions over the parameter domain, but construct them ``on-the-fly" during the evaluation stage, based on some optimality criteria. Another class of reduced order modelling techniques perform a spectral analysis of the training data, and use the part of the spectrum associated to high energy content to build the surrogate model. This is the case of the methods based on the Singular Value Decomposition (SVD) and their extensions (Proper Orthogonal Decomposition (POD) \cite{karhunen1947,loeve1963}, Balanced Truncation (see for instance \cite{antoulassorensen2001}), Multilinear Singular Value Decomposition, Proper Generalised Decomposition (PGD) \cite{ladevezepassieux2009,chinestaammarcueto2010,Nouy2010,chinestaammar2010,ladevezechamoin2011}), or on Krylov Subspaces (Moment Matching Method, see for instance \cite{antoulassorensen2001}). As in the case of interpolation-based reduced order models, the surrogate can be defined by an explicit spectral expansion over the parameter domain, like in the case of the POD. Otherwise, an implicit definition for the surrogate will require some ``online" computations to conform to an optimality criterion, like in the case of Galerkin-POD \cite{legresleyalonso2001,meyermatthies2003,kunishvolkwein2003,kunischvolkwein2011,kerfridengoury2012,gogupassieux2013} used in this article. In any case, the choice of the best reduced order modelling method for a particular application strongly depends on the characteristics of the underlying parametrised problem.

Reduced order modelling should be associated to error estimation to control the  distance between the solution of the parametrised problem (often called ``truth" solution) and the solution delivered by the surrogate. this requirement becomes even more stringent in the case of nested approximations, for instance in the case of multiscale modelling. It is fundamental to understand that there are two types of error associated with reduced order modelling: the ``offline" and the ``online" error. The ``offline" error is a global distance between the surrogate approximation and the ``truth" solution over the whole parameter domain. Error estimates for this quantity measure the average accuracy of the surrogate model over the parameter domain. This error can be evaluated by several means, depending on the application, which includes cross-validation estimates \cite{krzanowski1987,braconnierferrier2011,kerfridenschmidt2012,kerfridengoury2012}, estimates based on spectral analysis of the training data (see the review proposed in \cite{abdiwilliams2010}), or more advanced techniques that deliver bounds for a particular class of parametrised problems  \cite{rozzahuynh2008,ladevezechamoin2011}. These estimates are used to control the sampling of the parameter domain and the construction of the surrogate model. The ``online" error is the distance between the ``truth" solution and the solution delivered by the reduced model at a particular point of interest of the parameter domain. The accuracy and reliability of ``online" error estimates is crucial as they measure the quality of the solution that is actually delivered by the reduced model, which is in general not simply related to the average quality of the surrogate model \footnote{One of the fundamental ideas behind the Reduced Basis Method \cite{prudhommerovas2002} is to link the two types of errors by measuring the ``offline" error using a ``max"-type norm in the parameter domain. In general, the distinction between ``offline" and ``online" errors is not necessarily relevant if the actual output of the analysis is the quantity used as a target for the construction of the surrogate.}. In addition, the numerical complexity of the ``online" estimates should remain low for the evaluation stage to be performed at cheap costs. 

This paper focusses on the reliable, accurate and efficient bounding of the ``online" error in the context of the Galerkin-POD. In particular, we will focus on an elastostatic problem with discontinuous and parametrised elasticity constants, discretised by the Finite Element Method. The finite element mesh will be considered sufficiently fine so that the discretisation error can be neglected in comparison to the reduced order modelling error. The Snapshot-POD methodology will then be deployed to extract ``offline" an attractive spatial manifold, or reduced space, in which any solution to the parametrised problem of elasticity can be accurately represented. In the evaluation stage, an optimal solution corresponding to a particular set of elasticity constants can be optimally computed by a Galerkin projection of the governing equations in the reduced space. Upper bounding techniques for the reduced modelling error have been obtained for such problems in the context of the Reduced Basis Method. In \cite{rozzahuynh2008}, the error estimation relies on a Riesz representation for the parametrised residual, using a fixed bilinear form over the parameter domain. Then, the bounding property is obtained by weighing the result by a coercivity constant, which is a characteristic of the elliptic operator associated to the parametrised problem of interest. The evaluation of this constant is a key point of the strategy, and specific techniques such as the Successive Constraint Linear Optimization Method \cite{HuynhRozza2007} have been proposed to estimate this quantity whilst retaining the bounding property. However, coercivity lower bounds can be pessimistic in the case of elasticity. In this paper, we propose to proceed differently by using the concept of the Constitutive Relation Error (CRE) \cite{ladevezepelle2004}, which only requires to manipulate the concepts of displacement and stress admissibilities. In particular, no coercivity bound is required.

The CRE is a widely used technique to bound the error associated to a displacement-based approximation of elasticity problems. In particular, it has been applied to the evaluation of discretisation errors in a finite element context \cite{ladevezepelle1991}, where it coincides in practical implementations with the Equilibrated Residual approach (see for instance \cite{verfurth1988,ainsworthoden2000,steinruter2004,diezpares2010}).
Conceptually, the CRE proposes to construct a recovered stress field that is statically admissible, or equilibrated. Applying the constitutive relation to the kinematically admissible finite element solution that needs to be verified, one obtains a non-equilibrated stress field, called finite element stress field. The distance (in energy norm) between the recovered stress field and the finite element stress field is a bound for the discretisation error. All the technical difficulty resides in the construction of the equilibrated stress field \cite{pledchamoin2011}.

\begin{figure}[htb]
 \centering
 \includegraphics[width=1\linewidth]{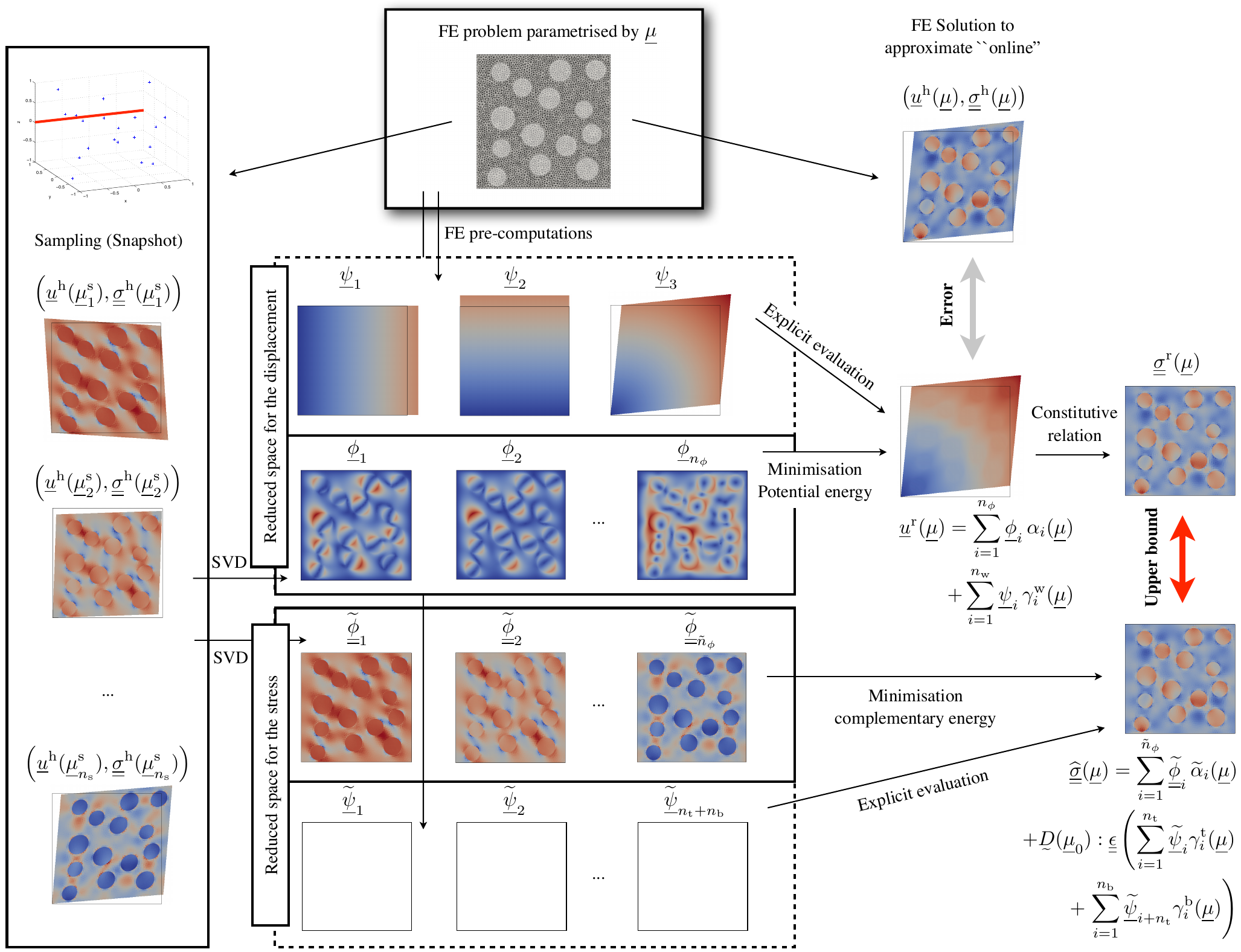}
 \caption{Schematic representation of the Galerkin-POD error bounding method based on the Constitutive Relation Error.}
 \label{fig:Schematic}
\end{figure}

We extend this idea to the certification of the Galerkin-POD, which will provide a bounding technique that is conceptually simple to understand, implement and control. Here, the reference is the finite element solution. Therefore, we first redefine the notion of statical admissibility, and require the recovered stress field to verify the equilibrium in the finite element sense. Then, at any point of the parameter domain, we can upper bound the reduced order modelling error by measuring a distance between this recovered field and stress field that is directly post-treated from the displacement delivered by the reduced model. In order for the recovered stress field to be available at cheap costs in the ``online" phase, we build a surrogate model for the finite element stress field. This reduced model is completely symmetric to the one developed for the displacement field (see figure \ref{fig:Schematic}, which can serve as visual guidance for the formal developments proposed in this paper). It requires an ``offline" training from the initial sampling of the parameter domain, and an ``online" computation to satisfy an optimality condition. While the optimality of the reduced displacement is enforced by minimisation of the potential energy, the optimality of the recovered stress is obtained by minimisation of the complementary energy in the space of stresses that are generated by the surrogate model. The technique proposed in this paper is largely influenced by, and  complementary to, the developments given in \cite{ladevezechamoin2011}, where the CRE is applied to evaluate the ``offline" error of a PGD reduced order modelling technique.

We will show that the efficiency of the upper bound can be seamlessly controlled by the order of the POD expansion of the reduced order model for the stress field. In order to lay the foundation for  the development of adaptive reduced order models, we also propose a lower bound, which requires the ``online" solution of a hierarchically enriched reduced model for the displacement field and whose effectivity can also be controlled. We will also provide some numerical results showing the convergence of both bounds with the refinement of the initial surrogate.

The efficiency of the method will be illustrated on problems that are directly related to computational homogenisation. However, the concepts presented in this paper apply to any affinely parametrised problems in linearised elasticity. There are two main reasons for focussing on homogenisation. Firstly, computational homogenisation can hugely benefit from algebra-based reduction methods, as mentioned earlier. A number of related contributions dealing with reduced order modelling for homogenisation acknowledge this fact \cite{yvonnethe2007,oskayfish2007,boyaval2008,abdullebai2013,ohlbergerschaefer2013}. Secondly, the homogenisation problem can, under some assumptions, be recast as a set of compliant boundary value problems over the RVE, which means that the energy in the domain is the actual output of the computation. Therefore, straight error estimates in the energy norm are of direct interest without resorting to specific goal-oriented techniques. This is not a limitation of the method itself, which could be coupled to any technique available in the literature to obtain general goal-oriented error estimates from error measures in the energy norm (see for instance \cite{odenprudhomme1999,larssonhansbo2002,steinruter2004,ladevezepelle2004,ladeveze2006,paresbonet2006,rozzahuynh2008}). This idea will be elaborated (but not used explicitly) in the core of this paper. Also, we limit our investigations to parametrised problems with fixed geometry (fixed distribution of inclusions in the context of homogenisation, and test sample that is large enough to be an RVE). Again this is not a limitation of the method, which could be coupled to the mapping technique presented in \cite{rozzahuynh2008} for instance.

The developments proposed in this contributions are organised as follows. In section  \ref{sec:Formulation}, we define the parametrised problem of elasticity and its discretisation by the Finite Element Method. In  section \ref{sec:MOR}, we present the basics of the Galerkin-POD reduced order modelling approach. The error bounding technique is developed in section \ref{sec:CRE}. Finally, a numerical validation of the certification methodology in the context of homogenisation will be presented in section \ref{sec:Results}.




\section{Parametric problem of elasticity solved by the Finite Element method}
\label{sec:Formulation}

\subsection{Problem statement: linear elasticity}

We formulate the problem of static equilibrium of a linear elastic structure occupying a bounded domain $\Domain$ in a physical space of dimension $\Dimension \in \{ 2, 3\}$. Let $\Point$ be an arbitrary point of domain $\Domain$ and let $\Pos= \GlobalCoordinate{1} \, \UnitVector{1} + \, ... \, + \GlobalCoordinate{d} \, \UnitVector{d}  $ be its cartesian decomposition in  reference frame $\ReferenceFrame=\left( \GlobalOrigin,  \, \UnitVector{1} , \, \UnitVector{2} , \, \UnitVector{3} \right)$. We look for a displacement field $\Disp \in \SpaceDisp = \Sobolev{1}(\Domain)$ that satisfies the Dirichlet boundary conditions $\Disp = \DirichletBC$ on the part $\DDomainDirichlet$ of the domain boundary $\DDomain$. Any displacement field that satisfies the conditions of regularity and the Dirichlet boundary conditions is said to be kinematically admissible and belongs to space $\SpaceDispAd \subset \SpaceDisp$. We introduce the Cauchy stress tensor field $\Stress$ which belongs to a space $\SpaceStress$ of sufficiently regular second-order tensor fields. A density of tractions $\NeumannBC$ is applied to the structure on the part $\DDomainNeumann \Def \DDomain \backslash \DDomainDirichlet$ of the domain. A density of forces denoted by $\VolumeForce$ is applied over $\Domain$. The principle of virtual work expresses the equilibrium of an arbitrary stress field belonging to $\SpaceStressAd \subset \SpaceStress$ as follows: 
\begin{equation}
\label{eqref:VirtualWork}
\forall \, \DispTest \in \SpaceDispZero, \qquad
- \int_\Domain \Stress : \StrainTest \, \DVolume
+
\int_{\Domain} \VolumeForce \cdot \DispTest  \, \DVolume +
\int_{\DDomainNeumann} \NeumannBC \cdot \DispTest   \, \DSurface  = 0 \, ,
\end{equation}
where $
\SpaceDispZero = \{ \DispAlt \in \SpaceDisp \, | \,  \DispAlt_{| \DDomainDirichlet} = \V{0}  \} $.
In the previous equation, $\Strain \Def \frac{1}{2} (\triangledown  \Disp + \triangledown  \Disp ^T)$ is the symmetric part of the displacement gradient. The solution to the problem of elasticity is an admissible pair $(\Disp,\Stress) \in \SpaceDispAd \times \SpaceStressAd $ that verifies the isotropic linear constitutive law
\begin{equation}
\label{eq:Consti}
\Stress = \LameFirst \, \text{Tr}( \Strain ) \, \IdentityTensor + 2 \, \LameSecond \, \Strain \Def \Hooke : \Strain \, ,
\end{equation}
where  $\LameFirst$ and $\LameSecond$ are the Lam\'e elasticity constants, and $\Hooke$ is the fourth-order Hooke's elasticity tensor. The inverse of this constitutive law reads
\begin{equation}
\Strain = \frac{1+\Poisson}{\Young} \, \text{Tr}( \Stress ) \, \IdentityTensor -\frac{\Poisson}{\Young} \, \Stress  \Def  \Compliance : \Stress \, ,
\end{equation}
with $\Young$ and $\Poisson$ the Young's and Poisson's moduli respectively, which are linked to the Lam\'e constants by the relationships $\LameFirst = \frac{\Young \, \Poisson}{(1+\Poisson)(1-2 \, \Poisson)}$ and $\LameSecond = \frac{\Young}{2(1+\Poisson)}$.

By substitution of the constitutive law into the principle of virtual work, the problem of elasticity can be recast in the primal variational form (displacement approach)
\begin{equation}
\displaystyle  \text{Find} \ \Disp \in \SpaceDispAd \ \text{such that} \ \forall \, \DispTest \in \SpaceDispZero,
\quad  
 a(\Disp, \DispTest) = l(\DispTest)\, ,
\end{equation}
where 
the symmetric bilinear form and the linear form associated with the problem of elasticity are respectively defined, for any fields $\DispAlt$ and $\DispTest$ in $\SpaceDisp$, by
\begin{equation}
\displaystyle 
a(\DispAlt,\DispTest) = \int_\Domain \StrainTest : \Hooke : \StrainBasic(\DispAlt)  \, \DVolume 
\, , \qquad 
l(\DispTest) =  \int_{\Domain} \VolumeForce \cdot \DispTest  \, \DVolume +
\int_{\DDomainNeumann} \NeumannBC \cdot \DispTest   \, \DSurface 
\, .
\end{equation}

\subsection{Parametrised problem of elasticity}

We consider that the input characterising the problem of elasticity are functions of a finite set of $\NumberParameter$ scalar parameters $(\ParameterI{i})_{i \in \llbracket 1,\NumberParameter\rrbracket}$ that are ordered in a parameter vector $\Parameter \in \mathbb{R}^{\NumberParameter}$. Let $\ParameterSpace \subset \mathbb{R}^{\NumberParameter}$ be the domain of admissibility of $\Parameter$. More precisely, 
the force densities $\VolumeForce$ and $\NeumannBC$, the Dirichlet boundary conditions $\DirichletBC$ and the Young's and Poisson's ratios  $\Young$ and $\Poisson$ are functions of parameter $\Parameter$ \footnote{The domain $\Domain$ and its boundary split $\DDomain = \DDomainDirichlet \cup \DDomainNeumann$ could also be parametrised, but this would lead to complications in terms of reduced order modelling and error estimation, which will not be addressed in this contribution. The interested reader can refer to the work presented in \cite{rozzahuynh2008} for the certification of Galerkin-based reduced order models applied to elliptic problems with parametrised geometries. In this work, specific mappings to a reference problem with fixed geometry are developed and used. This technique could be coupled in a relatively straightforward manner to the certified reduced order modelling methodology that we present in this paper}.

The field variables are now formally redefined as functions of the parameter:
\begin{equation}
\begin{array}{cccc}
\Disp : & \ParameterSpace & \rightarrow & \SpaceDisp
\\
 &  \Parameter & \mapsto & \DispParameter 
\end{array} 
\quad \text{and} \quad
\begin{array}{cccc}
\Stress : & \ParameterSpace & \rightarrow & \SpaceStress
\\
 &  \Parameter & \mapsto & \Stress(\Parameter) \, .
\end{array} 
\end{equation}
This definition is extended to all parametrised variables and spaces. If not stated explicitly, a symbolic expression of the type  $f(\Parameter) \in \mathcal{F}(\Parameter)$ will imply the definition of a function $f$ with inputs in $\ParameterSpace$ and values in space $\bigcup_{\Parameter \in \ParameterSpace} \mathcal{F}(\Parameter)$. A symbol of the type $f(y,z;\Parameter)$ will additionally imply that the values of $f$ are functions of variables $y$ and $z$.

For a given $\Parameter$, a kinematically admissible displacement field will be sought in $\SpaceDispAdParameter = \{  \DispAlt \in \SpaceDisp \, | \,  \DispAlt_{| \DDomainDirichlet} = \DirichletBCParameter \}$. Similarly, a statically admissible stress field $\Stress(\Parameter) \in \SpaceStressAdParameter$ will satisfy the parametrised principle of virtual work:
\begin{equation}
\label{eqref:VirtualWorkparametric}
\begin{array}{l}
 \displaystyle 
 \forall \, \DispTest \in \SpaceDispZero, \quad 
 \displaystyle
- \int_\Domain \Stress(\Parameter) : \StrainTest \, \DVolume
+
 \int_{\Domain} \VolumeForce(\Parameter) \cdot \DispTest  \, \DVolume +
 \int_{\DDomainNeumann} \NeumannBC(\Parameter) \cdot \DispTest   \, \DSurface = 0  \, ,
\end{array}
\end{equation}


By constraining the stress fields to satisfy \textit{a priori} the constitutive equation 
\begin{equation}
\label{eq:ConstiPara}
\Stress(\Parameter)= \Hooke(\Parameter) : \StrainBasic \left(\Disp(\Parameter) \right) \, , 
\end{equation}
the parametric problem of elasticity can be written in the primal variational form, for any $\Parameter \in \ParameterSpace$:
\begin{equation}
\label{eq:Bilineeee}
\begin{array}{l}
\displaystyle  \text{Find} \ \Disp(\Parameter) \in \SpaceDisp \ \text{such that}  \ \forall \, \DispTest \in \SpaceDispZero,
\\ \displaystyle  
 a(\Disp(\Parameter),\DispTest ;\Parameter) = l(\DispTest;\Parameter) \, ,
\end{array} 
\end{equation}
where the parametrised bilinear and linear forms associated with the problem of elasticity read, for any $\left(\DispAlt,\DispTest \right)$ in $\left( \SpaceDisp \right)^2$,
\begin{equation}
\displaystyle 
a (\DispAlt,\DispTest;\Parameter) = \int_{\Domain} \StrainBasic(\DispTest) : \Hooke(\Parameter) : \StrainBasic(\DispAlt) \, \DVolume
\, , \qquad  
l (\DispTest;\Parameter) =  \int_{\Domain} \VolumeForceParameter \cdot \DispTest  \, \DVolume +
\int_{\DDomainNeumann} \NeumannBCParameter \cdot \DispTest   \, \DSurface 
\, .
\end{equation}

\subsection{Quantity of interest and input-output map}
 
Let us consider a set of $n_\text{QoI}$ scalar quantities of interest ordered in a vector $\Output(\Parameter) = \widetilde{\Output}(\Disp(\Parameter)) \in \mathbb{R}^{n_\text{QoI}}$.  
Our goal is to construct a surrogate for map $\Output$.

The basic strategy applied in such cases is to compute the quantity of interest corresponding to certain (well-chosen) parameter values in $\Parameter \in \ParameterSpaceDiscr = \{ \Parameter^\text{s} _1 , \Parameter _2^\text{s} , \, ... \, , \Parameter^\text{s} _{n_\text{s}} \}$, where $\ParameterSpaceDiscr$ is a subset of $\ParameterSpace$ (called training parameter domain or snapshot parameter domain), and interpolate the quantity of interest in some ways over the entire parameter domain $\ParameterSpace$. Of course, the method of interpolation should be driven by considerations of optimality, stability and quality control. 

In particular, interpolating the solution $\Disp$ and not directly the output over the parameter domain is usually preferable. Notably, and this is the point that we focus on in this contribution, reliable error estimates can be obtained for a certain class of parametrised boundary value problems. We will use the Galerkin-POD approach, which will be described in section \ref{sec:MOR}.

\subsection{Finite element discretisation}

We approximate the solutions to the parametrised problem of elasticity by making use of a classical finite element discretisation $\SpaceDispFEM \subset \SpaceDisp$ of space $\SpaceDisp$ \cite{zienkiewicz1977,ciarlet1978}. As we assumed that domain $\Domain$ does not depend on parameter $\Parameter$, we can construct a unique finite element space for all the realisations of the parametrised boundary value problem. More precisely, the finite element space will be such that
\begin{equation}
\SpaceDispFEM = \left\{  \DispAlt \in \SpaceDisp \ | \ \forall \, j \in \{1, \, ... \, , \Dimension \}  , \ v_j \in \text{span} \left( (\ShapeFEMI{i})_{i \in \llbracket 1 , \NumberShape \rrbracket} \right) \right\} \, ,
\end{equation}
where $v_j$ denotes the $j^\text{th}$ component of vector field $\DispAlt$ and functions $(\ShapeFEMI{i})_{i \in \llbracket 1 , \NumberShape \rrbracket}$ are compactly supported finite element shape functions belonging to $\SpaceDisp$. 

Let $\SpaceDispFEMZero \Def \SpaceDispFEM \cap \SpaceDispZero $ be the space of finite element fields that vanish on $\DDomainDirichlet$ and let $\DispLift(\Parameter)$ be a particular field of $\SpaceDispAdParameter$, for any $\Parameter \in \ParameterSpace$. The finite element approximation $\DispFEM(\Parameter)$ of $\Disp(\Parameter)$ is the solution to the following variational problem:
\begin{equation}
\label{eq:FEMproblem1}
\begin{array}{l}
\displaystyle  \text{Find} \ \DispFEM(\Parameter)  \in \SpaceDispFEMZero + \{ \DispLift(\Parameter) \} \ \text{such that} \ \forall \, \DispTest \in \SpaceDispFEMZero,
\\ \displaystyle  
 a(\DispFEM(\Parameter),\DispTest;\Parameter) = l(\DispTest;\Parameter)   \, .
\end{array} 
\end{equation}

We will assume in this paper that the parametrised Dirichlet boundary conditions conform to the finite element space, which means that $\SpaceDispFEMAdParameter \Def \SpaceDispFEM \cap \SpaceDispAdParameter \neq \{ \}$. In this context, $ \DispLift(\Parameter) $ is naturally chosen in the finite element space $\SpaceDispFEM$ and we will use the alternative notation $\DispFEMLift(\Parameter) \equiv \DispLift(\Parameter)$, where $\DispFEMLift(\Parameter) \in \SpaceDispFEMAdParameter$. 

For any $\Parameter \in \ParameterSpace$, problem \eqref{eq:FEMproblem1} can finally be recast in the form:
\begin{equation} 
\label{eq:FEMproblem}
\begin{array}{l}
\displaystyle  \text{Find} \ \DispFEMHomo(\Parameter) \in \SpaceDispFEMZero \ \text{such that} \ \forall \, \DispTest \in \SpaceDispFEMZero,
\\ \displaystyle  
 a(\DispFEMHomo(\Parameter),\DispTest;\Parameter) = l(\DispTest;\Parameter) - a(\DispFEMLift(\Parameter),\DispTest;\Parameter)  \, ,
\end{array} 
\end{equation}
and the finite element solution is obtained making use of the lifting identity $\DispFEM(\Parameter) = \DispFEMHomo(\Parameter) + \DispFEMLift(\Parameter)$. Notice that the finite element solution obtained in this fashion is kinematically admissible, which is important for the remainder of the developments.


In the following, we assume that the finite element space is sufficiently fine so any measure the finite element error $\ErrorFEM(\Parameter) \Def  \Disp(\Parameter) - \DispFEM(\Parameter) $ is small for all $\Parameter \in \ParameterSpace$. 

\begin{CommentsFEM}
We introduce the following notations. Any field $\DispAlt \in \DispFEM$ will be represented by
\begin{equation}
\forall \, \Pos \in \Domain , \quad \DispFEMAlt(\Pos,\Parameter) = \MatrixShapeFEM(\Pos)[\DispFEMAlt(\Parameter)] \, ,
\end{equation}
where $\MatrixShapeFEM(\Pos) \in \mathbb{R}^\Dimension \times \mathbb{R}^\NumberUnknowns$ a matrix of shape functions, while $[\DispFEMAlt] \in \mathbb{R}^{\NumberUnknowns}$ is the vector of discrete unknowns corresponding to field $\DispFEMAlt$. We denoted by $\NumberUnknowns$ the number of discrete unknowns associated to the finite element discretisation. Similarly, the strain field $\StrainBasicEngin(\DispFEMAlt(\Parameter)$ corresponding to $\DispFEMAlt \in \SpaceDispFEM$, using the usual engineering notations, will be represented by
\begin{equation}
\forall \, \Pos \in \Domain , \quad \StrainBasicEngin(\DispFEMAlt(\Pos,\Parameter)) = \MatrixDerivativeShapeFEM(\Pos)[\DispFEMAlt(\Parameter)] \, ,
\end{equation}
where $\MatrixDerivativeShapeFEM(\Pos) \in  \mathbb{R}^\Dimension \in \mathbb{R}^{\widetilde{d}}$ is the matrix of shape function derivatives, where ${\widetilde{d}}$ is the number of components of the engineering strain vector field. Using these standard notations, problem \eqref{eq:FEMproblem} can be recast in the algebraic form, for any $\Parameter \in \ParameterSpace$
\begin{equation}
\label{eq:FEMproblemDiscr}
\begin{array}{l}
\displaystyle  \textit{Find} \ [\DispFEMHomo(\Parameter)] \in \bar{\mathcal{U}}^\text{0} \ \textit{such that} \ \forall \, [\DispTest] \in \bar{\mathcal{U}}^\text{0} \ \textit{that satisfies} \ \, ,
\\ \displaystyle  
[\DispTest]^T \StiffnessFEM(\Parameter) \, [\DispFEMHomo(\Parameter)] = [\DispTest]^T \left( \ForceFEM(\Parameter) - \StiffnessFEM(\Parameter) [\DispLift(\Parameter)]  \right)
\end{array} \, 
\end{equation}
where $\bar{\mathcal{U}}^\text{0} = \{  [\DispAlt] \in \mathbb{R}^\NumberUnknowns \  | \   \MatrixShapeFEM(\Pos) \, [\DispAlt] = \V{0}  , \, \forall \, \Pos \in  \DDomainDirichlet \}$ 
and the components of the finite element stiffness matrix and vector of external forces read:
\begin{equation}
\begin{array}{l}
\displaystyle
\StiffnessFEM(\Parameter)  = 
\int_\Domain \MatrixDerivativeShapeFEM^T \, \HookeEngin(\Parameter) \MatrixDerivativeShapeFEM \, \DVolume
\\ \displaystyle
\ForceFEM(\Parameter)  = \int_\Domain \MatrixShapeFEM^T \VolumeForce(\Parameter) \, \DVolume + \int_\DDomainNeumann \MatrixShapeFEM^T \NeumannBC(\Parameter) \, \DSurface
\end{array}
\end{equation}
In the previous equation, $\HookeEngin(\Parameter)$ is the Hooke matrix using engineering notations. The solution $\DispFEM(\Parameter)$ is post-processed as follows:
\begin{equation}
\forall \, \Pos \in \Domain , \quad \DispFEM(\Pos,\Parameter) = \MatrixShapeFEM(\Pos) \left(  [\DispFEMHomo(\Parameter)] + [\DispLift(\Parameter)]   \right)
\end{equation}
\end{CommentsFEM}

\section{Galerkin-Proper Orthogonal Decomposition}
\label{sec:MOR}



 
The idea of projection-based reduced order modelling relies on the fact that the solutions to parametrised boundary value problems are often found to lie in low-dimensional spatial subspaces engendered by the parametric dependence. In the training stage, several methods can be used to capture the attractive subspace numerically, amongst them the popular Snapshot (or Empirical) Proper Orthogonal Decomposition \cite{sirovich1987}.  
In the evaluation stage, the governing equations corresponding to any parameter of interest are projected in the reduced space obtained ``offline". In the case of elliptic parametrised boundary value problems, such a projection by means of the Galerkin method yields an optimally interpolated solution, from which the quantities of interest can be extracted.




\subsection{Projection-based reduced order modelling}


In order to introduce the approximation of the displacement $\DispFEM$ using a reduced space, let us first recall the lifting of the finite element solution over the parameter space
\begin{equation}
\label{eq:SplitDisplacement}
\forall \, \Parameter \in \ParameterSpace , \quad \DispFEM(\Parameter) = \DispFEMHomo(\Parameter) + \DispFEMLift(\Parameter) \, .
\end{equation}
In equation \eqref{eq:SplitDisplacement}, $\DispFEMHomo$ is a function of $\Parameter$ with values in $\SpaceDispFEMZero$, which means that the displacement fields $\DispFEMHomo(\Parameter)$ vanishes on $\DDomainDirichlet$ for any $\Parameter \in \ParameterSpace$, while the values of function $\DispFEMLift$ belong to the finite element space $\SpaceDispFEM$ and satisfy the Dirichlet boundary conditions of the parametric boundary value problem exactly. We assume for now that the lifting $\DispFEMLift$ is known, and concentrate on the approximation of the values of the homogeneous remainder $\DispFEMHomo$ using a reduced space.

Let us introduce a basis $( \ReducedBasisI{i} )_{ i \in \llbracket 1,  \NumberReducedBasis \rrbracket } \in \left( \SpaceDispFEMZero \right)^\NumberReducedBasis$ 
of the representative subspace $\SpaceDispRedZero \subset \SpaceDispFEMZero$ in which any value of the displacement $\DispFEMHomo$ will be approximated. This basis is also assumed to be known at this stage. Using these notations, and for any $\Parameter \in \ParameterSpace$, we look for an approximation $\DispReduced(\Parameter)$ of  $\DispFEM(\Parameter)$ to the parametrised problem of elasticity in the form
\begin{equation}
\label{eq:ReducedBasisApprox}
\displaystyle 
\DispFEM(\Parameter) \approx {\DispReduced}(\Parameter) \Def  \DispReducedHomo(\Parameter) + \DispFEMLift(\Parameter)  
\, , \quad 
\textrm{where} \quad \DispReducedHomo(\Parameter) = \sum_{i=1}^{\NumberReducedBasis} \ReducedBasisI{i} \, \ReducedCoeffI{i}(\Parameter) \, . 
\end{equation}
The functionals $(\ReducedCoeffI{i})_{i \in \llbracket 1,\NumberReducedBasis \rrbracket}$ are interpolation weights, called ``reduced variables''.
The interpolation weights can be optimally computed by using a Galerkin formulation of the elasticity problem in the reduced space, for any admissible parameter value. The projected problem corresponding to an arbitrary $\Parameter \in \ParameterSpace$ reads
\begin{equation}
\label{eq:GalerkinReducedBasis}
\begin{array}{l}
\displaystyle  \text{Find} \ \DispReducedHomo(\Parameter) \in \SpaceDispRedZero \ \text{such that} \ \forall \, \DispTest \in \SpaceDispRedZero,
\\ \displaystyle  
 a(\DispReducedHomo(\Parameter),\DispTest;\Parameter) = l(\DispTest;\Parameter) - a(\DispFEMLift(\Parameter),\DispTest;\Parameter) \, .
\end{array} 
\end{equation}

Substituting \eqref{eq:ReducedBasisApprox} into \eqref{eq:GalerkinReducedBasis}, and using the explicit form of the bilinear and linear forms of the parametrised elasticity problem, one obtains the following (small) system of $\NumberReducedBasis$ coupled equations in the reduced variables $(\ReducedCoeffI{i}(\Parameter))_{i \in \llbracket 1,\NumberReducedBasis \rrbracket}$:
\begin{equation}
\StiffnessFEM^\text{r}(\Parameter) \, \ReducedCoeff(\Parameter) = \ForceFEM^\text{r}(\Parameter) + \ForceFEM^\text{r,p}(\Parameter) \, ,
\end{equation}
where the $\ReducedCoeff(\Parameter) \in \mathbb{R}^\NumberReducedBasis$ is the unknown vector of weighting coefficients. The algebraic operators used in the previous expression are defined by:
\begin{equation}
\label{eq:SysReduced}
\begin{array}{ll}
\displaystyle
\forall \, (i,j) \in \llbracket 1 , \NumberReducedBasis \rrbracket^2 , 
& \displaystyle
 \StiffnessFEM_{ij}^\text{r}(\Parameter)   =  
a(\ReducedBasisI{j},\ReducedBasisI{i};\Parameter)
\\ \displaystyle
\forall \, j \in \llbracket 1 , \NumberReducedBasis \rrbracket ,
& \displaystyle
 \ForceFEM_{j}^\text{r}(\Parameter)  = 
l( \ReducedBasisI{j} ;\Parameter)
\\ \displaystyle
\forall \, j \in \llbracket 1 , \NumberReducedBasis \rrbracket , 
& \displaystyle
 \ForceFEM_j^\text{r,p}(\Parameter) =  
-a(\DispFEMLift(\Parameter),\ReducedBasisI{j};\Parameter) \, .
\end{array}
\end{equation}

\begin{CommentsFEM}
The algebraic operators can be linked to the finite element operators by substitution of the finite element interpolation for the reduced basis vectors $( \ReducedBasisI{i} )_{ i \in \llbracket 1,  \NumberReducedBasis \rrbracket }$:
\begin{equation}
\begin{array}{l}
\displaystyle 
\forall \, (i,j) \in \llbracket 1 , \NumberReducedBasis \rrbracket^2 , \quad \left. \StiffnessFEM^\text{r}(\Parameter) \right._{ij} = [\ReducedBasisI{j}]^T \StiffnessFEM(\Parameter) \, [\ReducedBasisI{i}]
\\ \displaystyle 
\forall \, j \in \llbracket 1 , \NumberReducedBasis \rrbracket , \quad
\left. \ForceFEM^\text{r}(\Parameter)  \right._j = [\ReducedBasisI{j}]^T  \ForceFEM(\Parameter)  
\\ \displaystyle 
\forall \, j \in \llbracket 1 , \NumberReducedBasis \rrbracket , \quad \left.  \ForceFEM^\text{p}(\Parameter) \right._j  = - [\ReducedBasisI{j}]^T  \StiffnessFEM(\Parameter) [\DispFEMLift(\Parameter)] \, .
\end{array}
\end{equation}
\end{CommentsFEM}

The linear system of equations \eqref{eq:SysReduced} can be solved inexpensively in the ``online'' phase, for any parameter value of interest. The result is the approximate field $\DispReduced(\Parameter) = \sum_{i=1}^{\NumberReducedBasis} \ReducedBasisI{i} \, \ReducedCoeffI{i}(\Parameter) + \DispFEMLift(\Parameter)$, from which the approximation of the output $\widetilde{\Output}(\DispReduced(\Parameter))$ can be estimated. However, in order to be able to compute each of the terms of system \eqref{eq:SysReduced} efficiently, some additional assumptions on the form of the parametrised boundary value problem are required, which is detailed in the next section. The computation of the reduced basis itself, and the definition of the lifting $\DispFEMLift$ is done in the ``offline'' phase, which will be detailed later on.

\subsection{Case of parametric elasticity problems admitting an affine form}

It is important to remark that a reduced order modelling technique cannot reach the expected numerical efficiency if the complexities of some of the operations performed in the ``online'' phase depend on the dimension of the underlying finite element space. The usual framework is to deal with problems that naturally admit an affine or radial form (see  \cite{rozzahuynh2008,ladevezechamoin2011} for instance). When the problem does not admit an affine form (e.g.: moving domains, nonlinear problems), an approximation of the initial governing equations must be performed in order to retrieve a reducible problem in affine form (see for instance \cite{barraultmaday2004,ryckelynck2008,niroomandialfaro2009,chaturantabutsorensen2010,carlbergbou-mosleh2011,kerfridengoury2012}). 

We will suppose in this contribution that the parametrised problem of interest \eqref{eq:FEMproblem} admits a natural affine form, which reads mathematically:
\begin{equation}
\label{eq:AffineForm}
\forall \, \Parameter \in \ParameterSpace , \, \forall \, (\DispTest, \DispAlt) \in (\SpaceDispFEMZero)^2, \quad
\begin{array}{l}
\displaystyle
a( \DispAlt, \DispTest ;\Parameter) = \sum_{i=1}^{\NumberBilinearAffine} \bar{a}_i ( \DispAlt, \DispTest ) \, \gamma_i^\text{a}(\Parameter) 
\\ \displaystyle
l( \DispTest ;\Parameter) = \sum_{i=1}^{\NumberLinearAffine} \bar{l}_i ( \DispTest) \, \gamma_i^\text{l}(\Parameter) 
\\ \displaystyle
a( \DispFEMLift(\Parameter) , \DispTest ;\Parameter) \Def l^\text{r,p}(\DispTest ;\Parameter) = \sum_{i=1}^{\NumberLinearLiftAffine} \bar{l}^\text{r,p}_i ( \DispTest) \, \gamma_i^\text{r,p}(\Parameter) \, ,
\end{array}
\end{equation}
where $\left( \bar{a}_i \right)_{i \in \llbracket 1, \NumberBilinearAffine \rrbracket}$ are parameter-independent bilinear forms (not necessarily symmetric positive definite), $\left( \gamma_i^\text{a}  \right)_{i \in \llbracket 1, \NumberBilinearAffine \rrbracket}$ are explicitly known functionals of the coefficients, 
$\left( \bar{l}_i     \right)_{i \in \llbracket 1, \NumberLinearAffine \rrbracket}$ are parameter-independent linear forms and $\left( \gamma_i^\text{l}  \right)_{i \in \llbracket 1, \NumberLinearAffine \rrbracket}$ are explicitly known functionals of the coefficients, 
$\left(  \bar{l}^\text{r,p}_i \right)_{i \in \llbracket 1, \NumberLinearLiftAffine \rrbracket}$ are parameter-independent linear forms and $\left(\gamma_i^\text{r,p}  \right)_{i \in \llbracket 1, \NumberLinearLiftAffine \rrbracket}$ are explicitly known functionals of the coefficients. Such a representation of the parametrised boundary value problem is obviously at hand if the parametrised data are originally given in the form of separate variables 
\begin{equation}
\label{eq:AffineQuantities}
\begin{array}{ll}
\displaystyle
\forall \, \Parameter \in \ParameterSpace, \,  \forall \, \Pos \in \Domain ,  \quad 
& \displaystyle
\Hooke(\Pos,\Parameter) = \sum_{i=1}^\NumberHookeAffine \HookeAffineI{i}(\Pos) \, \CoeffHookeAffineI{i}(\Parameter) 
\\ \displaystyle
\forall \, \Parameter \in \ParameterSpace, \, \forall \, \Pos \in \DDomainNeumann ,  \quad 
& \displaystyle
\NeumannBC(\Pos,\Parameter) = \sum_{i=1}^\NumberNeumannBCAffine \NeumannBCAffineI{i}(\Pos) \, \CoeffNeumannBCAffineI{i}(\Parameter) 
\\ \displaystyle
\forall \, \Parameter \in \ParameterSpace, \, \forall \, \Pos \in \Domain ,  \quad 
& \displaystyle
\VolumeForce(\Pos,\Parameter) = \sum_{i=1}^{\NumberVolumeForceAffine} \VolumeForceAffineI{i}(\Pos) \, \CoeffVolumeForceAffineI{i}(\Parameter) 
\\ \displaystyle
\forall \, \Parameter \in \ParameterSpace, \, \forall \, \Pos \in \DDomainDirichlet , \,  \quad 
& \displaystyle
\DirichletBC(\Pos,\Parameter) = \sum_{i=1}^{\NumberDirichletBCAffine} \DirichletBCAffineI{i}(\Pos) \, \CoeffDirichletBCAffineI{i}(\Parameter) \, ,
\end{array}
\end{equation}
where the notations used are consistent with those introduced for expressions \eqref{eq:AffineForm}.


In this context, assembling the projected system \eqref{eq:SysReduced} in the ``online'' phase can be done efficiently. The terms corresponding to each of the summands of the affine expansions \eqref{eq:AffineForm} can be pre-computed. For instance, assembling the projected stiffness $\StiffnessFEM^\text{r}(\Parameter)$, which is an operation of numerical complexity \textit{a priori} related to the dimension of the finite element space, can be written as follows:
\begin{equation}
\begin{array}{l}
\displaystyle 
\forall \, \Parameter \in \ParameterSpace, \quad \StiffnessFEM^\text{r}(\Parameter) = \sum_{k=1}^\NumberHookeAffine \bar{\StiffnessFEM}^\text{r}_k \, \CoeffHookeAffineI{k}(\Parameter)
\\ \displaystyle \text{where} , \quad 
 \forall \, k \in \llbracket 1 , \NumberHookeAffine \rrbracket , \, \forall \, (i,j) \in \llbracket 1 , \NumberReducedBasis \rrbracket^2 , 
 \,  \quad \quad \bar{\StiffnessFEM}^\text{r}_{k,ij} = 
\bar{a}_k( \ReducedBasisI{j} , \ReducedBasisI{j}) 
\Def \int_\Domain \StrainBasic(\ReducedBasisI{j}) : \HookeAffineI{k}  : \StrainBasic(\ReducedBasisI{i}) \, \DVolume  \, . 
\end{array}
\end{equation}
Operators $\left( \bar{\StiffnessFEM}^\text{r}_k \right)_{k \in \llbracket 1 , \NumberHookeAffine \rrbracket}$ can be precomputed ``offline''. In the ``online" phase, the assembly is simply done by computing a linear combination of these operators with coefficients $\left( \CoeffHookeAffineI{i}(\Parameter)  \right)_{i \in \llbracket 1, \NumberHookeAffine \rrbracket}$. The same technique can be used to assemble  $\ForceFEM^\text{r}(\Parameter)$ and $\ForceFEM^\text{r,p}(\Parameter)$, for any $\Parameter$ of interest. As a consequence, the numerical complexity of the ``online'' phase only depends on $\NumberHookeAffine$, $\NumberNeumannBCAffine$, $\NumberVolumeForceAffine$, $\NumberDirichletBCAffine$ and of course on the dimension $\NumberReducedBasis$ of the reduced space. 

\subsection{Non-homogeneous Dirichlet boundary conditions in projection-based reduced order modelling}

Applying non-homogeneous Dirichlet boundary conditions in projection-based reduced order modelling is not a trivial task. One possible approach is to use the usual boundary lifting performed in a finite element context. We propose an alternative approach based a global lifting that is (i) consistent with the choice of global basis vectors to define the reduced space $ \SpaceDispRedZero$ and (ii) consistent with the construction of the dual reduced space used for error estimation, as will be explained later on.

Our lifting technique makes use of the previously assumed form of the prescribed displacements \eqref{eq:AffineQuantities}. In the ``offline" stage, we compute a set of $\NumberDirichletBCAffine$ finite element displacement fields $(\FunctionAdKinematicAffineI{i})_{i \in \llbracket 1, \NumberDirichletBCAffine \rrbracket} \in \left(  \SpaceDispFEM \right)^\NumberDirichletBCAffine$ that satisfy
\begin{equation}
\label{eq:DispLift}
\forall \, i \in \llbracket 1, \, \NumberDirichletBCAffine   \rrbracket ,\quad 
\left\{ \begin{array}{l}
\displaystyle
\forall \, \DispTest \in \SpaceDispFEMZero , \quad a(\FunctionAdKinematicAffineI{i},\DispTest;\ParameterZero) = 0 
\\ \displaystyle
\forall \, \Pos \in \DDomainDirichlet , \quad
\FunctionAdKinematicAffineI{i}(\Pos ) = \DirichletBCAffineI{i}(\Pos) \, .
\end{array} \right.
\end{equation}
This set of fields is obtained by solving $\NumberDirichletBCAffine$ standard finite element problems ``offline'' with non-homogeneous boundary conditions. The choice of $\ParameterZero \in \ParameterSpace$ is arbitrary.

The finite element lifting function $\DispFEMLift$, which satisfies the Dirichlet boundary conditions for any $\Parameter \in \ParameterSpace$, can now be defined by the affine expansion
\begin{equation}
\label{eq:ParticularDisplacement}
\forall \, \Parameter \in \ParameterSpace , \quad \DispFEMLift(\Parameter) = \sum_{i=1}^{\NumberDirichletBCAffine} \FunctionAdKinematicAffineI{i} \, \CoeffDirichletBCAffineI{i}(\Parameter) \, .
\end{equation}

\subsection{Reduced spaces obtained by the Snapshot Proper Orthogonal Decomposition}
\label{sec:POD}


The purpose of the Snapshot Proper Orthogonal Decomposition \cite{sirovich1987} is to deliver an orthonormal reduced basis $( \ReducedBasisI{i} )_{ i \in \llbracket 1,  \NumberReducedBasis \rrbracket }$ of given cardinality $\NumberReducedBasis$ such that the $\mathcal{L}^2(\Domain)$-projection of a set of values of the ``truth" solution to the parametrised boundary value problem in the space generated by this basis is minimised. In our context, this optimisation problem reads:
\begin{equation}
\label{eq:SnapshotPOD}
\begin{array}{l}
\displaystyle   (\ReducedBasisI{i})_{ i \in \llbracket 1,  \NumberReducedBasis \rrbracket } 
= \underset{
( \V{\phi}_i^\star )_{ i \in \llbracket 1,  \NumberReducedBasis \rrbracket } \in \left( \SpaceDisp \right)^{\NumberReducedBasis} , \, \left<  \V{\phi}_i^\star , \V{\phi}_j^\star \right>_{\mathcal{L}^2(\Domain)}  = \, \delta_{ij} \ \forall \, (i,j) \in \llbracket 1 ,\NumberReducedBasis \rrbracket^2
}{\text{argmin}} \quad  J \left( ( \V{\phi}_i^\star )_{ i \in \llbracket 1,  \NumberReducedBasis \rrbracket } \right)
\\ \displaystyle
\text{where} \quad 
 J\left( ( \V{\phi}_i^\star )_{ i \in \llbracket 1,  \NumberReducedBasis \rrbracket } \right)
  = \frac{1}{\NumberReducedBasis} \sum_{\Parameter \in \ParameterSpaceDiscr} \left\| \DispFEMHomo(\Parameter) - \sum_{i=1}^\NumberReducedBasis \left< \DispFEMHomo(\Parameter) , \V{\phi}_i^\star \right>_{\mathcal{L}^2(\Domain)} \ \V{\phi}_i^\star \right\|_{\mathcal{L}^2(\Domain)}
\end{array} .
\end{equation}
In the previous equation, $< . \, , \, .  >_{\mathcal{L}^2(\Domain)} $ denotes the usual scalar product of $\mathcal{L}^2(\Domain)$, $\| \, . \,  \|_{\mathcal{L}^2}$ denotes the associated norm and $\delta_{ij}$ is the  Kronecker delta symbol. The snapshot parameter domain $ \ParameterSpaceDiscr = \{ \ParameterSnapshot{i} \, | \, i \in \llbracket 1, \NumberSnapshot \rrbracket \} \subset \ParameterSpace$ is a discrete set of $ \NumberSnapshot $ parameter values that are chosen in $\ParameterSpace$. The set of corresponding values of $\DispFEMHomo$ is the so-called snapshot $\SnapshotZero = \{  \DispFEMHomo(\Parameter) \, | \,  \Parameter \in \ParameterSpaceDiscr \}$, which is computed ``offline" using the ``truth" finite element approximation.

The solution to this classical problem of optimisation can be obtained by solving the eigenvalue problem
\begin{equation}
\CorrelationOperator \, \V{\zeta} = \lambda\,  \V{\zeta} 
\quad \text{where} \quad 
\bold{H}^\text{s}_{ij} = \left< {\DispFEMHomo}(\ParameterSnapshot{i}) , {\DispFEMHomo}(\ParameterSnapshot{j}) \right>_{\mathcal{L}^2(\Domain)} \  \ \forall \, (i,j) \in \llbracket 1 ,\NumberReducedBasis \rrbracket^2 \, .
\end{equation}
\begin{CommentsFEM}
Or in a discrete form by making use of the fact that the elements of the snapshot belong the the finite element space $\SpaceDispFEMZero$, and by calling the mass finite element matrix $\MassFEM = \int_\Domain \MatrixShapeFEM(\Pos)^T \MatrixShapeFEM(\Pos)  \, \DVolume$
\begin{equation}
\quad \bold{H}^\text{s}_{ij} = [\DispFEMHomo(\ParameterSnapshot{i})]^T \MassFEM \, [\DispFEMHomo(\ParameterSnapshot{j})] \, . 
\end{equation}
\end{CommentsFEM}
After ordering the eigenvalues $(\lambda_i)_{i \in \llbracket 1 , \NumberSnapshot \rrbracket}$ of $\CorrelationOperator$ in descending order and denoting by $(\V{\zeta}_i)_{i \in \llbracket 1 , \NumberSnapshot \rrbracket}$ the corresponding eigenvectors, the solution to the Snapshot POD optimisation problem \eqref{eq:SnapshotPOD} is given by
\begin{equation}
\forall \, i \in \llbracket 1 , \NumberReducedBasis \rrbracket, \quad \ReducedBasisI{i} = \sum_{j=1}^{\NumberSnapshot}  {\DispFEMHomo}(\ParameterSnapshot{j}) \,  \frac{\zeta_{i,j}}{\sqrt{\lambda_i}} \, ,
\end{equation}
\begin{CommentsFEM}
or, using the discrete notations,
\begin{equation}
\forall \, i \in \llbracket 1 , \NumberReducedBasis \rrbracket, \quad [\ReducedBasisI{i}] = \sum_{j=1}^{\NumberSnapshot}  [{\DispFEMHomo}(\ParameterSnapshot{j})] \,  \frac{\zeta_{i,j}}{\sqrt{\lambda_i}} \, ,
\end{equation}
\end{CommentsFEM}
where $\zeta_{i,j}$ denotes the $j^\text{th}$ component of eigenvector $\V{\zeta}_{i}$. For a given snapshot $\SnapshotZero$, the cardinality $\NumberReducedBasis$ of the reduced basis is usually chosen ``offline" such that the value of the objective function $ J\left( ( \ReducedBasisI{i} )_{ i \in \llbracket 1,  \NumberReducedBasis \rrbracket } \right)$ is small enough, which means that the correlated information contained in the snapshot has been correctly captured. Choosing the snapshot itself is a more difficult issue. Qualitatively, the sampling of the parameter domain should be fine enough to capture the variations of $\DispFEMHomo$ induced by the parametric dependence. In practice, this requirement can be validated by empirical means such as cross-validation  \cite{krzanowski1987,abdiwilliams2010,braconnierferrier2011,kerfridenschmidt2012}, or using more advanced numerical methods \cite{kunischvolkwein2011,tonurban2011,ladevezechamoin2011} to estimate some measure of the ``offline" error over the entire parameter domain.

\section{\textit{A posteriori} error bounding for projection-based reduced order modelling}
\label{sec:CRE}
 


\subsection{Error field and error measures}
\label{sec:QoI}
 
For any $\Parameter$ of interest, a relevant measure of the difference the finite element ``truth" solution $\DispFEM(\Parameter)$ and the displacement $\DispReduced(\Parameter)$ that is delivered by the Galerkin-POD needs to be evaluated. A practical measure of the error field $\ErrorReduced(\Parameter) \Def \DispFEM(\Parameter) - \DispReduced(\Parameter) $ can be obtained by evaluation of the error in the outputs $\widetilde{Q}_i(\ErrorReduced(\Parameter)) = \widetilde{Q}_i(\DispFEM(\Parameter)) - \widetilde{Q}_i(\DispReduced(\Parameter))$, for all $i \in \llbracket 1, n_\text{QoI} \rrbracket$. This error measure, usually referred to as error in the quantity of interest (QoI), can be related to an error measure in energy norm that can be bounded mathematically. Following classical work in this field \cite{odenprudhomme1999,steinruter2004,ladeveze2006,rozzahuynh2008,diezpares2010}, we can show\footnote{the proof is skipped as not central to this contribution} that the error in the QoI can be upper and lower bounded as follows: 
\begin{equation}
\label{eq:FrameQoI}
- \| \ErrorReduced(\Parameter)  \| _{\Hooke(\Parameter)} \| \V{e}_\text{z}^{\text{r},(i)}(\Parameter) \| _{\Hooke(\Parameter)} + R(\V{z}^{\text{r},(i)}(\Parameter);\Parameter) \leq \widetilde{Q}_i(\ErrorReduced(\Parameter)) \leq \| \ErrorReduced(\Parameter) \|_{\Hooke(\Parameter)} \| \V{e}_\text{z}^{\text{r},(i)}(\Parameter) \| _{\Hooke(\Parameter)}  + R(\V{z}^{\text{r},(i)}(\Parameter);\Parameter) \, .
\end{equation}
In the above inequalities, $\| \DispTest \|_{\Hooke(\Parameter)} = \left(\int_\Domain \StrainBasic(\DispTest ) : \Hooke(\Parameter) :  \StrainBasic(\DispTest) \, d \Domain \right)^\frac{1}{2}$ is the energy semi-norm associated to an arbitrary displacement field. $\V{e}_\text{z}^{\text{r},(i)}(\Parameter) \Def \V{z}^{\text{h},(i)}(\Parameter) - \V{z}^{\text{r},(i)} (\Parameter)$ is the error field of an auxiliary (or dual, adjoint) problem corresponding to quantity of interest $i$. The finite element ``truth" solution $\V{z}^{\text{h},(i)}(\Parameter) \in \SpaceDispFEMZero$ of this auxiliary problem  satisfies
\begin{equation}
\label{eq:DualProblem}
\forall \, \DispTest \in \SpaceDispFEMZero,
\quad 
 a(\DispTest,\V{z}^{\text{h},(i)}(\Parameter);\Parameter) =\widetilde{Q}_i(\DispTest) \, .
\end{equation}
This solution is as expensive to solve as the initial parametrised finite element problem. Its approximation $\V{z}^{\text{r},(i)}(\Parameter)$ in a reduced space $\mathcal{U}_\text{z}^{\text{r},0}(\Domain)$ spanned by finite element fields belonging to $\SpaceDispFEMZero$ satisfies
\begin{equation}
\forall \, \DispTest \in \mathcal{U}_\text{z}^{\text{r},0}(\Domain), \quad
 a(\DispTest,\V{z}^{\text{r},(i)}(\Parameter);\Parameter) = \widetilde{Q}_i(\DispTest) \, .
\end{equation}
Finally, the residual linear form $R$ of the initial ``truth" finite element problem and evaluated in the reduced order modelling solution is defined by
\begin{equation}
 \forall \, \DispTest \in \SpaceDisp, \quad
R( \DispTest ;\Parameter) \Def l(\DispTest;\Parameter) - a(\DispReduced(\Parameter)  , \DispTest ; \Parameter)  \, . 
\end{equation}
The residual term in \eqref{eq:FrameQoI}  is exactly computable in a reduced order modelling context (\textit{i.e.} the numerical complexity does not depend on the dimension of the finite element space). 
 
We will consider some important particular cases of the previous derivation: 
\begin{itemize}
\item if the reduced spaces used to approximate the initial and auxiliary problems are identical, \textit{i.e.} $\mathcal{U}_\text{z}^{\text{r},0}(\Domain) = \SpaceDispRedZero$, then the residual term in \eqref{eq:FrameQoI} vanishes owing to the Galerkin orthogonality property. In this case,  inequality set \eqref{eq:FrameQoI} reduces to the simple expression
\begin{equation}
| \widetilde{Q}_i(\ErrorReduced(\Parameter)) | \leq \| \ErrorReduced(\Parameter)  \| _{\Hooke(\Parameter)} \| \V{e}_\text{z}^{\text{r},(i)}(\Parameter) \| _{\Hooke(\Parameter)} \, ,
\end{equation}
\item if in addition to this, the linear form associated with quantity of interest $i$ is such that $\widetilde{Q}_i(\DispTest) = l(\DispTest;\Parameter) - a(\DispFEMLift(\Parameter),\DispTest;\Parameter) $, for all $\DispTest \in \SpaceDisp$ (\textit{i.e.} the initial problem \eqref{eq:FEMproblem} and the auxiliary problem \eqref{eq:DualProblem} are identical owing to the symmetry of bilinear form $a$), then we have
\begin{equation}
| \widetilde{Q}_i(\ErrorReduced(\Parameter)) | \leq  \| \ErrorReduced(\Parameter)  \|_{\Hooke(\Parameter)}^2 \, .
\end{equation}
The former case is often used when quantifying discretisation errors \cite{steinruter2004,gonzalesnadal2013} (the same mesh is used to solve the initial and auxiliary problems). More efficient bounding techniques have been proposed in this context, in particular in \cite{odenprudhomme1999}. The latter case does not require to solve an auxiliary problem. Although limited in terms of applications, we will see that this case arises naturally in the context of computational homogenisation.
\end{itemize}

In any case, the bounds defined by \eqref{eq:FrameQoI} are not computable as the energy norm of the exact error fields are not available. However, these quantities can be bounded from above an below. Providing that the bounds are sufficiently sharp and numerically affordable, they can be substituted in \eqref{eq:FrameQoI} to obtain guaranteed upper and lower bounds for the quantities of interest. We will employ the Constitutive Relation Error \cite{ladevezepelle2004} to bound the energy norm of the reduced order modelling error from above. This concept requires the availability of a stress field that satisfies the discrete principle of virtual work, or equilibrium in the finite element sense. Such a field can be obtained by ``offline" construction and ``online" evaluation of a POD-based surrogate for the finite element ``truth" stress over the parameter space. This surrogate is similar to the one developed for the approximation of the finite element ``truth" displacement, which makes the upper bounding technique conceptually simple. 

\subsection{Definition of the upper bound for the reduced order modelling error measured in energy norm}

We consider  a realisation of the parametrised problem of elasticity corresponding to an arbitrary parameter $\Parameter$ of $\ParameterSpace$, which is solved approximately using the reduced order modelling technique described in section \ref{sec:MOR}. The reduced order model delivers a kinematically admissible displacement field $\DispReduced(\Parameter) \in \SpaceDispFEMAd$. However, the stress field $\StressReduced(\Parameter) \Def \Hooke(\Parameter) :\StrainReduced(\Parameter)$ does not satisfy the equilibrium in the finite element sense \textit{a priori}. If it did, the solution to the ``truth" parametrised finite element problem $\DispFEM(\Parameter)$ would be at hand. The idea behind the proposed error bounding technique is to post-process a so-called ``recovered" stress field $\StressReducedHat(\Parameter) \in \SpaceStress$ that is equilibrated in the finite element sense. This admissibility condition reads:
\begin{equation}
\label{eq:StaticAdFEsense}
\forall \, \DispTest \in \SpaceDispFEMZero , \quad
- \int_\Domain \StressReducedHat(\Parameter) : \StrainTest \, \DVolume
+
 \int_{\Domain} \VolumeForce(\Parameter) \cdot \DispTest  \, \DVolume +
\int_{\DDomainNeumann} \NeumannBC(\Parameter) \cdot \DispTest   \, \DSurface = 0 \, .
\end{equation}
We denote by $\SpaceStressFEMAdParameter$ the space of stresses satisfying the parametrised equilbrium in the finite element sense \eqref{eq:StaticAdFEsense}.



\paragraph{Theorem:} 

The distance $\CRE(\Parameter)$ between the stress field $\StressReduced(\Parameter) \in \SpaceStress$ obtained by direct evaluation of the reduced order model and the recovered stress field $\StressReducedHat(\Parameter) \in \SpaceStressFEMAdParameter$ can be used to bound the error of reduced order modelling $\ErrorReduced(\Parameter) = \DispFEM(\Parameter) - \DispReduced(\Parameter) $ as follows:
\begin{equation}
 \CRE(\Parameter) \Def \| \StressReduced(\Parameter)  - \StressReducedHat(\Parameter)  \|_{\Compliance(\Parameter)}  \geq \| \ErrorReduced(\Parameter) \|_{\Hooke(\Parameter)}  \, ,
\end{equation}
where $ \| \StressAlt \|_{\Compliance(\Parameter)} = \left(  \int_\Domain \StressAlt : \Compliance(\Parameter) : \StressAlt \, \DVolume  \right)^\frac{1}{2}$ is the energy norm associated to an arbitrary stress field $\StressAlt \in \SpaceStress$.

\paragraph{Proof:}

The proof is a straightforward extension of the one used to bound the discretisation error in the context of finite element approximations \cite{ladevezepelle2004}. We start by expanding the distance between the reduced stress field $\StressReduced(\Parameter)$ and the recovered stress field $\StressReducedHat(\Parameter)$ by using the trivial identity
\begin{equation}
\label{eq:Proof1}
\| \StressReduced(\Parameter) - \StressReducedHat(\Parameter)   \|^2_{\Compliance(\Parameter)} = \left\| \left( \StressReduced(\Parameter) - \StressFEM(\Parameter) \right)+ \left( \StressFEM(\Parameter)  - \StressReducedHat(\Parameter) \right) \right\|^2_{\Compliance(\Parameter)} \, ,
\end{equation}
where $\StressFEM(\Parameter) \Def \Hooke(\Parameter) :\StrainBasic\left( \DispFEM(\Parameter) \right)$ is the ``truth" finite element stress field (\textit{i.e.}: the stress field which would be obtained without reduced order modelling). Then, by using the constitutive relation and the definition of the energy norms given previously, we can expand identity \eqref{eq:Proof1} as follows:
\begin{equation}
\label{eq:CREDeveloped}
\begin{array}{ll}
\displaystyle
\| \StressReduced(\Parameter) - \StressReducedHat(\Parameter)  \|^2_{\Compliance(\Parameter)} = 
& \displaystyle
  \| \DispReduced(\Parameter) - \DispFEM(\Parameter)  \|_{\Hooke(\Parameter)}^2 + 
  \| \StressFEM(\Parameter)  - \StressReducedHat(\Parameter)   \|_{\Compliance(\Parameter)}^2 
\\
& \displaystyle
+ 2 \int_\Domain
\left( \StressFEM(\Parameter)  - \StressReducedHat(\Parameter) \right) : \left( \StrainBasic \left(\DispReduced(\Parameter) \right) - \StrainBasic \left( \DispFEM(\Parameter) \right)  \right)
 \, \DVolume  \, .
 \end{array}
\end{equation}
Recall that both the finite element stress field and the recovered stress field are equilibrated in the finite element sense (equation \eqref{eq:StaticAdFEsense}). Therefore, the following identity holds:
\begin{equation}
\forall \, \DispTest \in \SpaceDispFEMZero, \quad 
\int_\Domain \,  \left( \StressFEM(\Parameter) -  \StressReducedHat(\Parameter) \right) : \StrainTest  \DVolume = 0 \, .
\end{equation}
As both the finite element displacement field and the displacement field obtained by solving the reduced order model belong to the finite element space $\SpaceDispFEM$ and satisfy the Dirichlet boundary conditions, we obtain that the last term in \eqref{eq:CREDeveloped} vanishes, which yields
\begin{equation}
\label{eq:CREWithEfficiencyTerm}
\left( \CRE(\Parameter) \right)^2 =  
  \|  \ErrorReduced(\Parameter)  \|_{\Hooke(\Parameter)}^2 
  +  \| \StressFEM(\Parameter)  -  \StressReducedHat(\Parameter) \|_{\Compliance(\Parameter)}^2 \, .
\end{equation}
The last term of the right-hand side of the above equation is positive, which concludes the proof. \hfill $\square$ \\

A first important remark is that the Galerkin orthogonality property $a( \ErrorReduced(\Parameter) , \DispTest; \Parameter) = 0$ for all $\DispTest \in \SpaceDispRedZero$ is not used, which means that the bounding property  would hold true if an interpolation method other than Galerkin reduced order modelling had been used to interpolate $\DispReduced$ over the parameter domain (\textit{e.g.} Kriging, interpolation over \textit{a priori} defined polynomial bases, \textit{etc.}). However, $\DispReduced(\Parameter)$ would still be required to be kinematically admissible.

A second important remark is that the efficiency of the error estimate depends on the distance $\|  \StressFEM(\Parameter)  -  \StressReducedHat(\Parameter)  \|_{\Compliance(\Parameter)}$,(\textit{i.e.}: the distance between the recovered stress field and the ``truth" finite element stress field), for any $\Parameter \in \ParameterSpace$, which needs to be kept in mind when constructing the equilibrated stress field. If the recovered stress field is finite element stress field, the CRE estimate $\CRE(\Parameter)$ coincides with the exact error. Of course, the finite element stress field is not affordable in the ``online" phase. Therefore, the next question is \textit{''how can we compute a recovered stress field that is
\begin{itemize}
\item a good approximation of the finite element stress field,
\item equilibrated in the finite element sense,
\item of ``online'' numerical complexity approximately equal to the numerical complexity of evaluating the reduced order model?''
\end{itemize} 
}

\subsection{Principle of the construction for the equilibrated stress fields}

In order to address the previous question, we propose to build a POD-based reduced order model for the finite element stress. The technique is similar to the one used to build the reduced order model for the displacement field. First, we use the affine representation of the external load to build a statically admissible stress field over the whole parameter domain. Then, the complementary part of the finite element stress field defined over the parameter domain is sampled, and a basis for the subspace spanned by these samples is extracted using a singular value decomposition. These operations are performed ``offline''. In the ``online'' stage, for a given parameter, the coefficients associated to the reduced basis functions are optimally obtained by solving a projected problem.




Specifically, we start by splitting the finite element stress into two parts, as follows
\begin{equation}
\label{eq:SplitRecoveredStress}
\forall \, \Parameter \in \ParameterSpace, \quad \StressFEM(\Parameter) = \StressFEMHomo(\Parameter) +  \StressFEMLift(\Parameter) \, .
\end{equation}
The first part $\StressFEMHomo$ belongs to a space $\SpaceStressFEMZero$ of stress fields satisfying the homogeneous equilibrium conditions associated with the ``truth" finite element problem, which reads
\begin{equation}
\label{eq:EquiliZero}
\forall \, \Parameter \in \ParameterSpace, \, \forall \, \DispTest \in \SpaceDispFEMZero, \quad 
\int_\Domain \StressFEMHomo(\Parameter) : \StrainTest \, \DVolume = 0 \, .
\end{equation}
The second part of split \eqref{eq:SplitRecoveredStress} is a particular stress field $\StressFEMLift(\Parameter) \in \SpaceStressFEMAd$ that satisfies the equilibrium in the finite element sense (equation \eqref{eq:StaticAdFEsense}).
$\StressFEMLift$ will be explicitly defined as a function of the parameter, while the complementary part $\StressFEMHomo$ will be approximated using the Snapshot-POD:
\begin{equation}
\forall \, \Parameter \in \ParameterSpace, \quad \StressFEM(\Parameter)  \approx \StressReducedHat(\Parameter) \Def  \StressFEMReducedHomo(\Parameter) +  \StressFEMLift(\Parameter) \, ,
\end{equation}
where the approximate stress $\StressFEMReducedHomo(\Parameter) \in \SpaceStressFEMZero$ is such that it satisfies the homogeneous equilibrium equations in the finite element sense for any $\Parameter$ in $\ParameterSpace$.

It is clear that within this framework, the recovered stress field $\StressReducedHat(\Parameter)$ satisfies the equilibrium in the finite element sense \eqref{eq:StaticAdFEsense} over the entire parameter domain.


\subsection{Riesz representation for the parametrised static load}

Let us first construct the particular stress $\StressFEMLift$ associated to non-homogeneous equilibration conditions. We make use of the assumed affine form of the Neumann boundary conditions and body forces given in expression \eqref{eq:AffineQuantities}. In the ``offline" phase, we compute a set of global finite element vectors $(\FunctionAdStaticAffineI{i})_{i \in \llbracket 1, \, \NumberNeumannBCAffine + \NumberVolumeForceAffine   \rrbracket} \in \left( \SpaceDispFEMZero \right)^{\NumberNeumannBCAffine + \NumberVolumeForceAffine}$ corresponding to the summand of the affine form identities by solving successively the finite element problems:
\begin{equation}
\label{eq:Riesz}
\begin{array}{ll}
\displaystyle
\forall \, i \in \llbracket 1, \, \NumberNeumannBCAffine   \rrbracket ,\, \forall \, \DispTest \in \SpaceDispFEMZero , 
& \displaystyle
a(\FunctionAdStaticAffineI{i},\DispTest;\ParameterZero) = 
\int_\DDomainNeumann \NeumannBCAffineI{i} \cdot \DispTest \, \DSurface
\\ \displaystyle
\forall \, i \in \llbracket 1, \,  \NumberVolumeForceAffine    \rrbracket ,\, \forall \, \DispTest \in \SpaceDispFEMZero , 
& \displaystyle
a(\FunctionAdStaticAffineI{i+\NumberNeumannBCAffine},\DispTest;\ParameterZero) = \int_\Domain \VolumeForceAffineI{i} \cdot \DispTest \, \DVolume 
\end{array}
\end{equation}
Notice that we have chosen to enforce homogeneous Dirichlet conditions for this series of finite element problems, which ensures a certain regularity of the procedure (the parametrised bilinear form $a( \, . \, , \, . \,  ; \Parameter_0)$ is positive definite over $\left( \SpaceDispFEMZero \right)^2$). Parameter $\ParameterZero$ is the one that we also used to define the lifting $\DispFEMLift$. This choice has been made for the sake of simplicity.


In the ``online'' phase of the error estimation procedure, the fields $(\FunctionAdStaticAffineI{i})_{i \in \llbracket 1, \, \NumberNeumannBCAffine + \NumberVolumeForceAffine   \rrbracket}$
 are avilable, and for a given parameter $\Parameter \in \ParameterSpace$, we evaluate the field
  $\RieszLift(\Parameter) \in \SpaceDispFEMZero$ defined by the formula
\begin{equation}
\RieszLift(\Parameter) =
 \sum_{i=1}^\NumberNeumannBCAffine \FunctionAdStaticAffineI{i} \,  \CoeffNeumannBCAffineI{i}(\Parameter) + \sum_{i=1}^{\NumberVolumeForceAffine} \FunctionAdStaticAffineI{i+\NumberNeumannBCAffine} \,  \CoeffVolumeForceAffineI{i}(\Parameter) \, .
\end{equation}
We can then verify that the stress field 
\begin{equation}
\label{eq:StressParticular}
\StressFEMLift(\Parameter) \Def \Hooke(\ParameterZero) : \StrainBasic(\RieszLift(\Parameter)) \, ,
\end{equation}
is statically admissible in the finite element sense, for any $ \Parameter \in \ParameterSpace $.
Therefore, the quantity $\| \StressFEMLift(\Parameter) - \StressReduced(\Parameter) \|_{\Compliance(\Parameter)}$ is an upper bound for the error measure $\| \ErrorReduced(\Parameter)  \|_{\Hooke(\Parameter)}$. However, in the case where the body load and applied tractions are zero, this bound is trivial of no practical interest. In the general case anyway, this bound should be sharpened by computing the complement $\StressFEMReducedHomo(\Parameter)$ to the recovered stress, which is done next by making use use of an additional spectral analysis of the training data.

\subsection{Snapshot POD for the finite element stress}

The ``offline'' computation of the snapshot delivers a set of admissible stress field in the finite element sense $\SnapshotStress \Def \{ \StressFEM(\Parameter) \, | \, \Parameter \in \ParameterSpaceDiscr   \}$. After subtracting the corresponding values of the previously defined statically admissible stress component $\StressFEMLift$, we obtain the set  $\SnapshotStressZero \Def  \{ \StressFEMHomo(\Parameter) \, | \, \Parameter \in \ParameterSpaceDiscr  \}$ of stress fields satisfying the homogeneous equilibrium equations in the finite element sense. In the ``online'' phase, our purpose is to compute the correction $\StressFEMReducedHomo(\Parameter)$ corresponding to an arbitrary parameter $\Parameter \in \ParameterSpace$ as an optimal combination of these fields in the sense of the minimisation of  $\| \StressReducedHat(\Parameter) - \StressFEM(\Parameter) \|_{\Compliance(\Parameter)} = \| \StressFEMReducedHomo(\Parameter) - \StressFEMHomo(\Parameter) \|_{\Compliance(\Parameter)}$. In other words, we aim at maximising the efficiency of the error estimate.

In order to do so, we compute a singular value decomposition of the finite element stress fields contained in the snapshot $\SnapshotStressZero$. We look for a set of $\NumberReducedBasisStress$ basis tensor fields $( \ReducedBasisStressI{i} )_{ i \in \llbracket 1,  \NumberReducedBasisStress \rrbracket } \in (\SpaceStress)^\NumberReducedBasisStress$ that are orthogonal with respect to the inner product $< . \, ,  . >_{\Compliance(\ParameterZero)} = \int_\Domain .   : \Compliance(\ParameterZero) :  . \, \DVolume$ of space $\SpaceStress$, and that are solution to the following optimisation problem:

\begin{equation}
\label{eq:SnapshotPODStress}
\begin{array}{l}
\displaystyle   (\ReducedBasisStressI{i})_{ i \in \llbracket 1,  \NumberReducedBasisStress \rrbracket } 
= \underset{
( \M{\widetilde{\phi}}_i^\star )_{ i \in \llbracket 1,  \NumberReducedBasis \rrbracket } \in \left( \SpaceStress \right)^{\NumberReducedBasisStress} , \, < \M{\widetilde{\phi}}_i^\star , \, \M{\widetilde{\phi}}_j^\star  >_{\Compliance(\ParameterZero)}= \delta_{ij} \ \forall \, (i,j) \in \llbracket 1 ,\NumberReducedBasisStress \rrbracket^2
}{\text{argmin}} \quad  \widetilde{J} \left( (\M{\widetilde{\phi}}_i^\star )_{ i \in \llbracket 1,  \NumberReducedBasisStress \rrbracket } \right)
\\ \displaystyle
\text{where} \quad 
\widetilde{J} \left( ( \M{\widetilde{\phi}}_i^\star  )_{ i \in \llbracket 1,  \NumberReducedBasisStress \rrbracket } \right)
  = \frac{1}{\NumberReducedBasisStress} \sum_{\Parameter \in \ParameterSpaceDiscr} 
  \left\| \StressFEMHomo(\Parameter) - \sum_{i=1}^\NumberReducedBasisStress \left< \StressFEMHomo(\Parameter) , \M{\widetilde{\phi}}_i^\star \right>_{\Compliance(\ParameterZero)} \ \M{\widetilde{\phi}}_i^\star \right\|_{\Compliance(\ParameterZero)}
\end{array} .
\end{equation}

This problem is similar to problem \eqref{eq:SnapshotPOD} except that we deal with second order tensors, and that the optimality of the decomposition is defined in the sense of the weighted norm $\| \, . \,  \|_{\Compliance(\ParameterZero)}$ associated to $< . \, ,  . >_{\Compliance(\ParameterZero)}$ instead of the usual Frobenius norm. The energy norm used to define the optimality of the POD is evaluated in $\ParameterZero$ to avoid adding unnecessary parameters to the bounding technique. Similarly, we have assumed that the sampling of the parameter domain performed to train the surrogates for the stress and for the displacement are the same in order to keep the methodology simple. 

The solution to optimisation problem \eqref{eq:SnapshotPODStress} is obtained by solving the eigenvalue problem
\begin{equation}
\label{eq:CorrelationOpStress}
\CorrelationOperatorStress \, \widetilde{\V{\zeta} } = \widetilde{\lambda} \,  \widetilde{\V{\zeta} }
\quad \text{where} \quad 
\bold{\widetilde{H}}^\text{s}_{ij} = \left< {\StressFEMHomo}(\ParameterSnapshot{i}) , {\StressFEMHomo}(\ParameterSnapshot{j}) \right>_{\Compliance(\ParameterZero)} \  \ \forall \, (i,j) \in \llbracket 1 ,\NumberReducedBasisStress \rrbracket^2 \, .
\end{equation}
After arranging the eigenvalues $(\widetilde{\lambda}_i)_{i \in \llbracket 1 , \NumberSnapshot \rrbracket}$ of $\CorrelationOperatorStress$ in descending order and denoting by $(\V{\widetilde{\zeta}}_i)_{i \in \llbracket 1 , \NumberSnapshot \rrbracket}$ the corresponding eigenvectors, the solution to the Snapshot POD optimisation problem \eqref{eq:SnapshotPODStress} is given by
\begin{equation}
\forall \, i \in \llbracket 1 , \NumberReducedBasisStress \rrbracket, \quad \ReducedBasisStressI{i} = \sum_{j=1}^{\NumberSnapshot}  {\StressFEMHomo}(\ParameterSnapshot{j}) \,  \frac{\widetilde{\zeta}_{i,j}}{\sqrt{\widetilde{\lambda}_i}} \, .
\end{equation}






Finally, the recovered stress field is given over the whole parameter domain by the surrogate model
\begin{equation}
\label{eq:recovered}
\forall \, \Parameter \in \ParameterSpace , \quad \StressReducedHat(\Parameter)  =  \StressFEMReducedHomo(\Parameter) +  \StressFEMLift(\Parameter) 
= \sum_{i=1}^\NumberReducedBasisStress \ReducedBasisStressI{i} \, \ReducedCoeffStressI{i}(\Parameter)  +  \StressFEMLift(\Parameter) \, ,
\end{equation}
where  $(\ReducedCoeffStressI{i})_{i \in \llbracket 1 , \NumberReducedBasisStress \rrbracket} \in \mathbb{R}^{\NumberReducedBasisStress}$ are interpolation functionals, and are the only unknowns to be estimated in the ``online'' phase. 

Notice that $\StressFEMReducedHomo(\Parameter)$ satisfies the homogeneous equilibrium equations for any $\Parameter \in \ParameterSpace$ because the reduced basis stress fields $( \ReducedBasisStressI{i} )_{ i \in \llbracket 1,  \NumberReducedBasisStress \rrbracket }$ are linear combinations of fields that satisfy the homogeneous equilibrium equations. Then, by linearity, the recovered stress $\StressReducedHat(\Parameter)$ is indeed equilibrated in the finite element sense.

\subsection{``Online" evaluation of the reduced order model for the stress}

In the ``online'' phase, a simple method to estimate the coefficients $(\ReducedCoeffStressI{i}(\Parameter))_{i \in \llbracket 1 , \NumberReducedBasisStress \rrbracket}$ of the recovered stress evaluated at a particular parameter value $\Parameter \in \ParameterSpace$ is to perform an \textit{a priori} interpolation of these coefficients over the parameter domain. This could be done, for instance, by means of a moving least-squares technique or Kriging. 

A more advanced interpolation technique consists in computing a recovered stress $\StressReducedHat(\Parameter)$ optimally in the sense of the maximisation of the efficiency of the error estimate $\CRE(\Parameter)$. In other words, for a particular $\Parameter \in \ParameterSpace$, we look for a recovered stress field that is compatible with the surrogate model and is solution to the optimisation problem:
\begin{equation}
\StressReducedHat(\Parameter) = \underset{ \StressAlt \in \SpaceStressReducedParameter  }{\text{argmin}} \|  \StressAlt - \StressFEM(\Parameter) \|_{\Compliance(\Parameter)} \, ,
\end{equation}
where the space of admissibility for the reduced stress field corresponding to parameter $\Parameter$ is defined by $\SpaceStressReducedParameter = \{ \StressAlt \in \SpaceStress  \, | \, \StressAlt = \sum_{i=1}^\NumberReducedBasisStress \ReducedBasisStressI{i} \, \alpha^\star_i  +  \StressFEMLift(\Parameter) , \, \forall \, (\alpha^\star_i)_{i \in \llbracket 1,\NumberReducedBasisStress \rrbracket } \in \mathbb{R}^{\NumberReducedBasisStress} \}$. 
This problem of optimisation can be recast in the variational form
\begin{equation}
\label{eq:MiniStress}
\begin{array}{l}
\displaystyle
\text{Find} \ \StressReducedHat(\Parameter) \in \SpaceStressReducedParameter \ \text{such that} \ \forall \, \StressTest \in \SpaceStressReducedZero, 
\\ \displaystyle
- \int_\Domain \StressReducedHat(\Parameter) : \Compliance(\Parameter) :  \StressTest  \, \DVolume +  
\int_\Domain \StrainBasic \left(\DispFEM(\Parameter) \right) :  \StressTest  \, \DVolume = 0 \, ,
\end{array}
\end{equation}
where the constitutive relation has been used to obtain the right-hand side of the equation, and the space $\SpaceStressReducedZero$ is defined by $\SpaceStressReducedZero = \{ \StressAlt \in \SpaceStress  \, | \, \StressAlt = \sum_{i=1}^\NumberReducedBasisStress \ReducedBasisStressI{i} \, \alpha^\star_i  , \, \forall \, (\alpha^\star_i)_{i \in \llbracket 1,\NumberReducedBasisStress \rrbracket } \in \mathbb{R}^{\NumberReducedBasisStress} \}$.

Now, by recalling that $\StressReducedHat(\Parameter) = \StressFEMReducedHomo(\Parameter) +  \StressFEMLift(\Parameter) $ and that $\DispFEM(\Parameter) = \DispFEMHomo(\Parameter) +  \DispFEMLift(\Parameter)$ and taking into account equation \eqref{eq:EquiliZero} and the fact that $\DispFEMHomo(\Parameter) \in \SpaceDispFEMZero$, we obtain the following variational form for the determination of the recovered stress field: 
\begin{equation}
 \forall \, \StressTest \in \SpaceStressReducedZero,  \quad 
\int_\Domain \StressFEMReducedHomo(\Parameter) : \Compliance(\Parameter) :  \StressTest  \, \DVolume = \int_\Domain  \StrainBasic \left(\DispFEMLift(\Parameter) \right) :  \StressTest  \, \DVolume
- \int_\Domain \StressFEMLift(\Parameter) : \Compliance(\Parameter) :  \StressTest  \,  \DVolume
 \, .
\end{equation}

Therefore, the expansion coefficients $(\ReducedCoeffStressI{i}(\Parameter))_{i \in \llbracket 1 , \NumberReducedBasisStress \rrbracket}$ of the reduced stress $\StressFEMReducedHomo(\Parameter)$ are obtained by solving the linear system of equations
\begin{equation}
\label{eq:SysReducedStress}
\widetilde{\StiffnessFEM}^\text{r}(\Parameter) \, \widetilde{\ReducedCoeff}(\Parameter) = \widetilde{\ForceFEM}^\text{r}(\Parameter) + \widetilde{\ForceFEM}^\text{r,p}(\Parameter)  \, ,
\end{equation}
where the components of vector $\widetilde{\ReducedCoeff}(\Parameter)$ are the interpolation coefficients $\left( \ReducedCoeffStressI{i}(\Parameter) \right)_{i \in \llbracket 1 , \NumberReducedBasisStress \rrbracket}$, and 
\begin{equation}
\label{eq:SysReducedStress2}
\begin{array}{ll}
\displaystyle
\forall \, (i,j) \in \llbracket 1 , \NumberReducedBasisStress \rrbracket^2 , 
& \displaystyle
 \widetilde{\StiffnessFEM}_{ij}^\text{r}(\Parameter)   =  \int_\Domain \ReducedBasisStressI{j} : \Compliance(\Parameter)  : \ReducedBasisStressI{i} \, \DVolume 
\\ \displaystyle
\forall \, j \in \llbracket 1 , \NumberReducedBasisStress \rrbracket , 
& \displaystyle
\widetilde{\ForceFEM}_j ^\text{r}(\Parameter) = \int_\Domain  \ReducedBasisStressI{j}  : \StrainBasic \left(\DispFEMLift(\Parameter) \right)   \, \DVolume 
\\ \displaystyle
\forall \, j \in \llbracket 1 , \NumberReducedBasisStress \rrbracket , 
& \displaystyle
\widetilde{\ForceFEM}_j ^\text{r,p}(\Parameter) = - \int_\Domain \ReducedBasisStressI{j}  : \Compliance(\Parameter)  :  \StressFEMLift(\Parameter)     \, \DVolume 
 \, . 
\end{array}
\end{equation}

\paragraph{``Offline"/``Online" split of the numerical complexity.}

Assembling linear system \eqref{eq:SysReducedStress} can be done efficiently by a combination of ``offline" and ``online" computations, where all the operations with numerical complexity that depends on the dimension of the underlying finite element space are performed ``offline". 

We will first assume that the parametrised compliance tensor $\Compliance$ admits the natural affine form:
\begin{equation}
\forall \, \Parameter \in \ParameterSpace, \, \forall \, \Pos \in \Domain, \quad \Compliance(\Pos ; \Parameter) = \sum_{i=1}^{\NumberComplianceAffine} \ComplianceAffineI{i}(\Pos) \, \CoeffComplianceAffineI{i}(\Parameter) \, .
\end{equation}
Then, the expression of operator $\widetilde{\StiffnessFEM}^\text{r}$ can be expanded as follows: 
\begin{equation}
\label{eq:AssemblyStress1}
\begin{array}{l}
\displaystyle \forall \, \Parameter \in \ParameterSpace, \quad 
\widetilde{\StiffnessFEM}^\text{r}(\Parameter) = \sum_{k=1}^\NumberComplianceAffine \bar{\widetilde{\StiffnessFEM}}^\text{r}_k \, \CoeffComplianceAffineI{k}(\Parameter)
\\ \displaystyle \text{where} , \quad
\forall \, k \in \llbracket 1 , \NumberComplianceAffine \rrbracket , \, \forall \, (i,j) \in \llbracket 1 , \NumberReducedBasisStress \rrbracket^2 , \quad\bar{\widetilde{\StiffnessFEM}}^\text{r}_{k,ij} =   \int_\Domain \ReducedBasisStressI{j} : \ComplianceAffineI{k}  : \ReducedBasisStressI{i} \, \DVolume  \, .
\end{array}
\end{equation}
By making use of the expansion of quantity $\DispFEMLift$ over the parameter domain (equation \eqref{eq:ParticularDisplacement}, we find that the first term in the right hand side of linear system \eqref{eq:SysReducedStress} reads 
\begin{equation}
\label{eq:AssemblyStress2}
\begin{array}{l}
\displaystyle \forall \, \Parameter \in \ParameterSpace, \quad 
\widetilde{\ForceFEM}^\text{r}(\Parameter) =  \sum_{k=1}^\NumberDirichletBCAffine  \bar{\widetilde{\ForceFEM}}^\text{r}_{k} \, \CoeffDirichletBCAffineI{k}(\Parameter)
\\ \displaystyle \text{where} , \quad
\forall \, k \in \llbracket 1 , \NumberDirichletBCAffine \rrbracket , \, \forall \, j \in \llbracket 1 , \NumberReducedBasisStress \rrbracket , \quad \bar{\widetilde{\ForceFEM}}^\text{r}_{k,j} =   \int_\Domain \StrainBasic(\FunctionAdKinematicAffineI{k}) : \ReducedBasisStressI{j} \, \DVolume  \, .
\end{array}
\end{equation}
The same approach, using \eqref{eq:StressParticular}, permits to obtain the affine expansion of $\widetilde{\ForceFEM} ^\text{r,p}(\Parameter)$ over the parameter domain.

The``online" assembly operations, comprising the first line of \eqref{eq:AssemblyStress1} and \eqref{eq:AssemblyStress2} respectively, do not depend on the dimension of the finite element space. The remaining operations, namely the integrations specified by the second line of \eqref{eq:AssemblyStress1} and \eqref{eq:AssemblyStress2} respectively, are performed ``offline". 

\subsection{Computation of the upper bound}

Once coefficients $(\ReducedCoeffStressI{i})_{i \in \llbracket 1 , \NumberReducedBasisStress \rrbracket}$ corresponding to a particular parameter $\Parameter \in \ParameterSpace$ have been computed, the ``online" error estimate $\CRE(\Parameter)$ defined by
\begin{equation}
\left( \CRE(\Parameter) \right)^2 = \int_\Domain \left( \StressReduced(\Parameter) - \StressReducedHat(\Parameter) \right) : \Compliance(\Parameter) : \left( \StressReduced(\Parameter) - \StressReducedHat(\Parameter) \right)   \, \DVolume \, ,
\end{equation}
can be evaluated. This can be done by a combination of ``online" and ``offline" computations that make use of the affine form of $\Compliance$, the affine expansions of $\StressFEMReducedHomo$ (equation \eqref{eq:recovered}) and $\StressFEMLift$ (equation \eqref{eq:StressParticular}) over the parameter domain, and the affine expansions of $\DispReducedHomo$ (equation \eqref{eq:ReducedBasisApprox}) and $\DispFEMLift$ (equation \eqref{eq:ParticularDisplacement}). The technique is similar to the one deployed to assemble system \eqref{eq:SysReducedStress} and will not be detailed for the sake of concision.

\subsection{Lower bound}

As will be shown in the results section of this paper, deriving an efficient lower bound for  the ``online" error associated to the reduced order modelling strategy can help control the efficiency of the upper bounding method. In particular, the choice of the dimension $\NumberReducedBasisStress$ of the reduced space for the recovered stress remains to be made at this stage.

A lower bound $\LowerBound(\Parameter)$ for the error measure $\| \ErrorReduced(\Parameter) \|_{\Hooke(\Parameter)}$ corresponding to any parameter $\Parameter \in \ParameterSpace$ can be obtained by first constructing an enhanced surrogate for the displacement:
\begin{equation}
\forall \Parameter \in \ParameterSpace, \quad 
\DispFEM(\Parameter) \approx \DispReducedEnhanced(\Parameter) = \DispFEMLift(\Parameter) + \DispReducedEnhancedHomo(\Parameter) \, , 
\end{equation}
such that $\DispReducedEnhancedHomo(\Parameter)$ belongs to a reduced space $\SpaceDispReducedEnhancedZero$ richer than $\SpaceDispRedZero$ (\textit{i.e.} $\SpaceDispRedZero \subset \SpaceDispReducedEnhancedZero \subset \SpaceDispFEMZero$). In order to define this enriched space, we simply perform ``offline" a Snapshot-POD of higher order for the displacement field, which reads
\begin{equation}
\SpaceDispReducedEnhancedZero \Def \text{span} \left( (\ReducedBasisI{i} )_{i \in \llbracket 1, \NumberReducedBasisEnriched \rrbracket } \right) \, ,
\end{equation}
where $\NumberSnapshot \geq \NumberReducedBasisEnriched > \NumberReducedBasis$ and, consistently with the notations of section \ref{sec:POD}, 
\begin{equation}
\forall \, i \in \llbracket \NumberReducedBasis+1 , \NumberReducedBasisEnriched \rrbracket, \quad \ReducedBasisI{i} = \sum_{j=1}^{\NumberSnapshot}  {\DispFEMHomo}(\ParameterSnapshot{j}) \,  \frac{\zeta_{i,j}}{\sqrt{\lambda_i}} \, .
\end{equation}

In the ``online" phase, $\DispReducedEnhanced(\Parameter)$ is optimally obtained by a Galerkin projection of the governing equations in space $\SpaceDispReducedEnhancedZero$, which reads
\begin{equation}
\label{eq:GalerkinReducedBasisEnriched}
\begin{array}{l}
\displaystyle  \text{Find} \ \DispReducedEnhancedHomo(\Parameter) \in \SpaceDispReducedEnhancedZero \ \text{such that} \ \forall \, \DispTest \in \SpaceDispReducedEnhancedZero,
\\ \displaystyle  
 a(\DispReducedEnhancedHomo(\Parameter),\DispTest;\Parameter) = l(\DispTest;\Parameter) - a(\DispFEMLift(\Parameter),\DispTest;\Parameter) \, .
\end{array} 
\end{equation}
The field $\DispReducedEnhanced(\Parameter)$ is a hierarchically enhanced approximation of $\DispFEM(\Parameter)$. It is expected to be closer than $\DispReduced(\Parameter)$, in some sense, to the ``truth" finite element solution $\DispFEM(\Parameter)$. The method to assemble and solve problem \eqref{eq:GalerkinReducedBasisEnriched} using a combination of ``offline" and ``online" operations is the same as the one detailed in section \ref{sec:MOR}, the only difference being the dimensions of the respective reduced spaces.

\paragraph{Theorem}
Defining an approximation of the error $\ErrorReduced$ by $\ErrorReducedApprox(\Parameter) \Def   \DispReducedEnhanced(\Parameter) - \DispReduced(\Parameter)$ for any $\Parameter \in \ParameterSpace$, we have the following lower bounding property
\begin{equation}
\forall \,  \Parameter \in \ParameterSpace , \quad \LowerBound(\Parameter) \Def \frac{ \displaystyle
|  R( \ErrorReducedApprox(\Parameter) ; \Parameter)  |} 
{  \displaystyle
\|  \ErrorReducedApprox(\Parameter)  \|_{\Hooke(\Parameter)}
} 
\leq 
\| \ErrorReduced(\Parameter)  \|_{\Hooke(\Parameter)} \, ,
\end{equation}
where the residual $R$ is the parametrised linear form defined by 
\begin{equation} 
\forall \,  \Parameter \in \ParameterSpace , \, \forall \, \DispTest \in \SpaceDisp, \quad
R( \DispTest ;\Parameter) \Def l(\DispTest;\Parameter) - a(\DispReduced(\Parameter)  , \DispTest ; \Parameter)  \, . 
\end{equation}
\paragraph{Proof} The proof is an extension of a similar technique used in the finite element context (see for instance \cite{diezpares2010}). The starting point is the weak form \eqref{eq:FEMproblem} of the parametrised problem of elasticity, from which we can obtain, for any $\Parameter \in \ParameterSpace$,
\begin{equation}
\forall \, \DispTest \in \SpaceDispFEMZero, \quad 
a(\DispFEMHomo(\Parameter)-\DispReducedHomo(\Parameter),\DispTest;\Parameter) = l(\DispTest ;\Parameter ) - a(\DispFEMLift(\Parameter),\DispTest ;\Parameter) - a(\DispReducedHomo(\Parameter),\DispTest;\Parameter) \, .
\end{equation} 
Using the identities $\DispReduced = \DispReducedHomo + \DispFEMLift $ and $\DispFEM = \DispFEMHomo + \DispFEMLift  $, we obtain the weak form of the equations governing the error  $\ErrorReduced$:
\begin{equation}
\forall \, \DispTest \in \SpaceDispFEMZero, \quad 
a(\ErrorReduced(\Parameter),\DispTest;\Parameter) = R(\DispTest ;\Parameter ) \, .
\end{equation}
We can now substitute for $\DispTest$ the approximate error $\ErrorReducedApprox(\Parameter) \in \SpaceDispFEMZero$, which leads to the expression
\begin{equation}
\forall \, \DispTest \in \SpaceDispFEMZero, \quad 
a(\ErrorReduced(\Parameter),\ErrorReducedApprox(\Parameter);\Parameter) =R( \ErrorReducedApprox(\Parameter) ;\Parameter ) \, .
\end{equation}
The exact error $\ErrorReduced(\Parameter)$ and its approximation $\ErrorReducedApprox(\Parameter)$ belong to the finite element space $\SpaceDispFEMZero$. The bilinear form $a$ is an inner product for this particular space. We can therefore apply the Cauchy-Schwartz inequality to obtain
\begin{equation}
\sqrt{a(\ErrorReduced(\Parameter),\ErrorReduced(\Parameter);\Parameter)} \sqrt{a(\ErrorReducedApprox(\Parameter),\ErrorReducedApprox(\Parameter)} \geq | R( \ErrorReducedApprox(\Parameter) ;\Parameter ) | \, 
\end{equation}
and the announced result is immediate by making use of the definition of $\| \, . \,  \|_{\Hooke(\Parameter)}$. \hfill $\square$

Again, an efficient combination of ``offline" and  ``online" computations can be deployed such that the ``online" numerical complexity associated with the evaluation of $\LowerBound$ does not depend on the dimension of the finite element space. 


Summarising the upper and lower bounding results developed in the previous sections, we have, for any $\Parameter \in \ParameterSpace$
\begin{equation}
\LowerBound(\Parameter) \leq \| \ErrorReduced(\Parameter) \|_{\Hooke(\Parameter)} \leq \CRE(\Parameter) \, .
\end{equation}


\section{Example: homogenisation of random composite materials}
\label{sec:Results}

\begin{figure}[htb]
 \centering
 \includegraphics[width=0.68\linewidth]{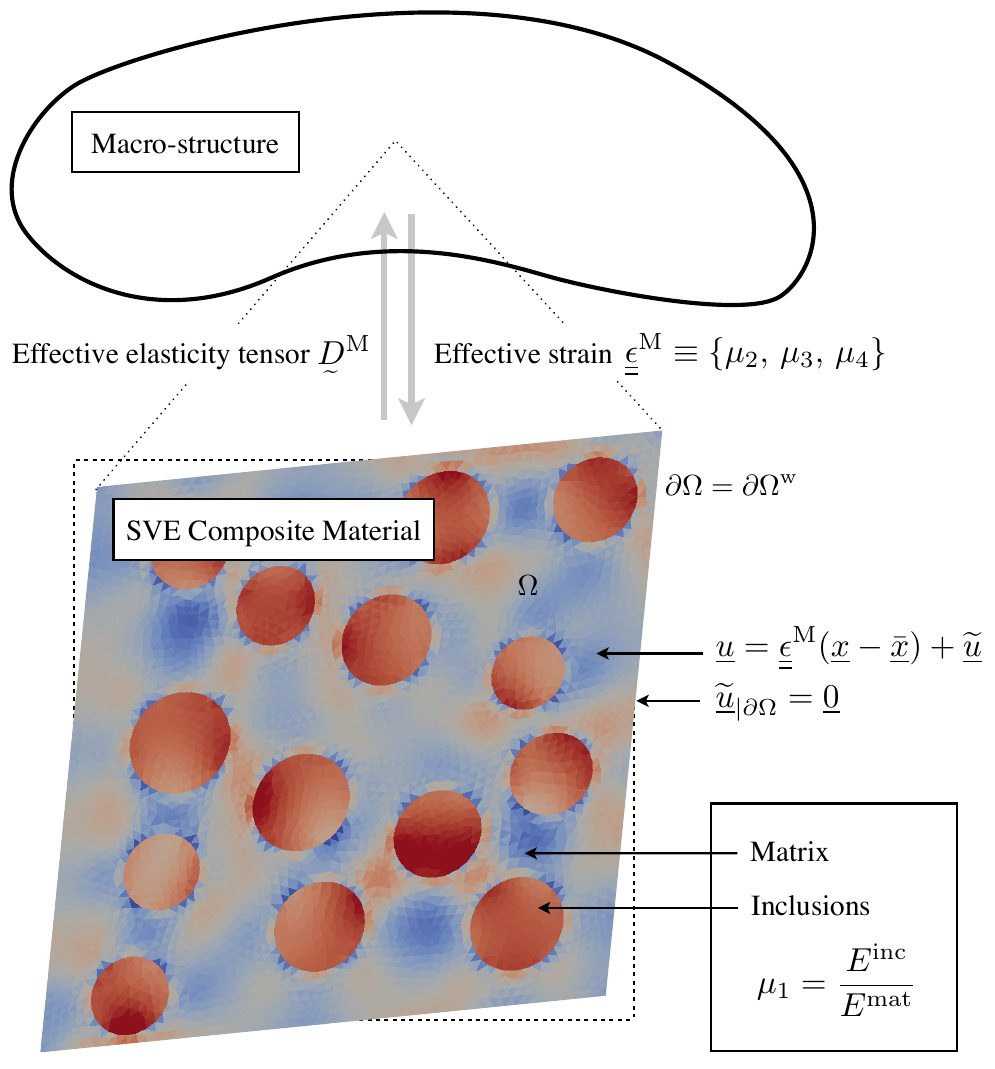}
 \caption{Schematic representation of the computational homogenisation framework for composite materials. The elastic constrast of the particulate composite is parametrised.}
 \label{fig:Homo}
\end{figure}

In order to illustrate the efficiency of the certified reduced order modelling methodology described previously, we will apply it in the context of a computational homogenisation scheme for random composite materials, in dimension 2 (plane strain assumption). The heterogeneous material of interest is made of two isotropic, linear elastic phases possessing distinct elastic constants: circular inclusions and surrounding matrix (see figure \ref{fig:Homo}). The positions and diameters of the inclusions are distributed randomly. The aim is to determine 
the so-called effective elasticity tensor as a function of some characteristics $\Parameter^\text{m}  \in \ParameterSpace^\text{m}$ of the material heterogeneities. In other words, we want to build a virtual chart of the homogenised properties of the class of composite materials under investigation.


\subsection{Classical computational homogenisation scheme for composite materials}
 
Following classical approaches in homogenisation of random heterogeneous media (see for instance \cite{nemat-nasser1999}), we consider that domain $\Domain$ is a statistical volume element (SVE) of the composite material. The SVE is loaded via uniform Dirichlet boundary conditions. Parameters $\Parameter^\text{l} \in \ParameterSpace^\text{l}$ fully characterise this load. We define the complete set of parmameters $\Parameter = \left( \left. \Parameter^\text{m} \right.^T \ \left. \Parameter^\text{l} \right.^T \right)^T$. For a given material characterised by $\Parameter^\text{m}$, the displacement field is additively split into a ``coarse'' contribution $\bar{\Disp}$ and a fluctuation field $\widetilde{\Disp}$:
\begin{equation}
\label{eq:SplitHomo}
\forall \, \Parameter \in \ParameterSpace, \, \forall \, \Pos \in \Domain , \quad \Disp(\Pos;\Parameter) = \bar{\Disp}(\Pos;\Parameter) + \widetilde{\Disp}(\Pos;\Parameter) \, . 
\end{equation}
The ``coarse'' contribution $\bar{\Disp}$ lies in a low-dimensional ``smooth" space $\SpaceDispCoarse$ of dimension 3 defined by
\begin{equation}
\SpaceDispCoarse = \{ \DispTest \in \SpaceDisp \, | \, \DispTest(\Pos) = \StrainBasic^\text{M}  (\Pos - \bar{\Pos}) , \, \forall \, \StrainBasic^\text{M} \in \mathbb{R}^2 \times \mathbb{R}^2 \ \text{such that} \ \StrainBasic^\text{M} =\left. \StrainBasic^\text{M} \right.^T  \} \, ,
\end{equation}
where  $\bar{\Pos}$ is the barycenter of $\Domain$ (\textit{i.e.} $\int_\Domain (\Pos - \bar{\Pos}) \,  \DVolume = \V{0}$). The effective strain tensor $\StrainBasic^\text{M}$ represents a far-field action applied to the material. The complementary fluctuation field $\widetilde{\Disp}$ lies in space $\SpaceDispFluctuation$ of fields with small characteristic length of variation defined by:
\begin{equation}
\label{eq:SpaceFluctu}
\SpaceDispFluctuation = 
 \left\{ \DispTest \in \SpaceDisp \, | \, \left< \StrainBasic(\DispTest) \right>_\Domain=0   \right\} \, ,
\end{equation}
where the averaging operator is defined by $\left<  \ . \ \right>_\Domain \Def \frac{1}{|\Domain |} \int_\Domain  \ . \  \DVolume$. In order the displacement field to be completely determined under the action of the coarse field, we enforce that the fluctuation field vanishes on the boundary $\DDomain \equiv \DDomainDirichlet$ of the SVE, consistently with definition \eqref{eq:SpaceFluctu} and the so-called macrohomogeneity condition (see \cite{nemat-nasser1999} for reference):
\begin{equation}
\label{eq:BCSVE}
\forall \, \Parameter \in \ParameterSpace, \, \forall \, \Pos \in \DDomain, \quad {\Disp}(\Pos;\Parameter) = \bar{\Disp}(\Pos;\Parameter)  \, .
\end{equation}

Assuming zero body force, the parametrised boundary value problem \{ \eqref{eqref:VirtualWorkparametric},
\eqref{eq:ConstiPara}, \eqref{eq:BCSVE} \} associated with the equilibrium of the SVE is well-posed and for any effective strain. The effective elasticity tensor $\Hooke^\text{M}$ can now be defined by the relationship
\begin{equation}
\begin{array}{l}
\label{eq:EffectiveDef}
\displaystyle
\forall (\StressTest,\DispTest) \in \SpaceStressAd  \times \left( \SpaceDispCoarse + \SpaceDispFluctuation \right) \ \text{satisfying} \
\eqref{eq:BCSVE} \ \text{and} \ \eqref{eq:ConstiPara} , 
\\ \displaystyle
\left<  \StressTest \right>_\Domain = \Hooke^\text{M}(\Parameter^\text{m}) \, \left< \StrainBasic(\DispTest) \right>_\Domain \, . 
\end{array}
\end{equation}

The components of the effective elasticity tensor can be obtained by numerical testing (see for instance \cite{zohdiwriggers2005}), whereby one applies a range of elementary effective strains through the Dirichlet boundary conditions, then solves the corresponding SVE boundary value problem numerically and finally post-processes the resulting stress. More precisely, individual components of the effective Hooke tensor are obtained by prescribing effective strain $\StrainBasic^\text{M}(\Parameter^\text{l}) = \frac{1}{2} \left( \UnitVector{k} \otimes \UnitVector{l} +  \UnitVector{l} \otimes \UnitVector{k} \right)$, for $(k,l) \in \{ 0,1 \}^2$, solve the associated Dirichlet boundary value problem for displacement field $\Disp(\Parameter)$, and finally use the extractor $\M{\Sigma} \Def \frac{1}{2} \left( \UnitVector{i} \otimes \UnitVector{j} +  \UnitVector{j} \otimes \UnitVector{i} \right)$ for $(i,j) \in \{ 0,1 \}^2$ to compute
\begin{equation}
D^\text{M}_{ijkl}(\Parameter^\text{m}) = \M{\Sigma} : \left<  \Hooke(\Parameter^\text{m}) : \StrainBasic(\Disp(\Parameter)) \right>_\Domain \, .
\end{equation}
In our case, the circular inclusions are distributed in an isotropic manner. For a sufficiently large SVE, the effective elastic law will be linear isotropic (equation \eqref{eq:Consti} applied to the average stress and strain tensors). The effective Lam\'e constants $\lambda^\text{M}(\Parameter^\text{m})$ and $G^\text{M}(\Parameter^\text{m})$ can be extracted by performing two numerical tests and computing the quantity of interest
\begin{equation}
Q(\Parameter) = \widetilde{Q}(\Disp(\Parameter)) = \StrainBasic^\text{M}(\Parameter^\text{l}) : \left<  \Hooke(\Parameter^\text{m}) : \StrainBasic(\Disp(\Parameter)) \right>_\Domain \, .
\end{equation}
Indeed, choosing $ \StrainBasic^\text{M}(\Parameter^\text{l}) = \frac{1}{2} \left( \UnitVector{1} \otimes \UnitVector{2} +  \UnitVector{2} \otimes \UnitVector{1} \right)$ and following this procedure yields $Q(\Parameter) = G^\text{M}(\Parameter^\text{m})$ directly, while choosing $ \StrainBasic^\text{M}(\Parameter^\text{l}) = \UnitVector{1} \otimes \UnitVector{1} $ gives us $Q(\Parameter) = \lambda^\text{M}(\Parameter^\text{m}) + 2 \, G^{M}(\Parameter^\text{m})$. 

Still, solving the SVE boundary value problem for arbitrary parameters of the heterogeneities can require a tremendous numerical effort. We reduce this effort by deploying the Galerkin-POD described in section \ref{sec:MOR} for the parametrised SVE problem.

\subsection{Data and discretisation of the parametrised SVE problem}

\begin{figure}[htb]
 \centering
 \includegraphics[width=0.58\linewidth]{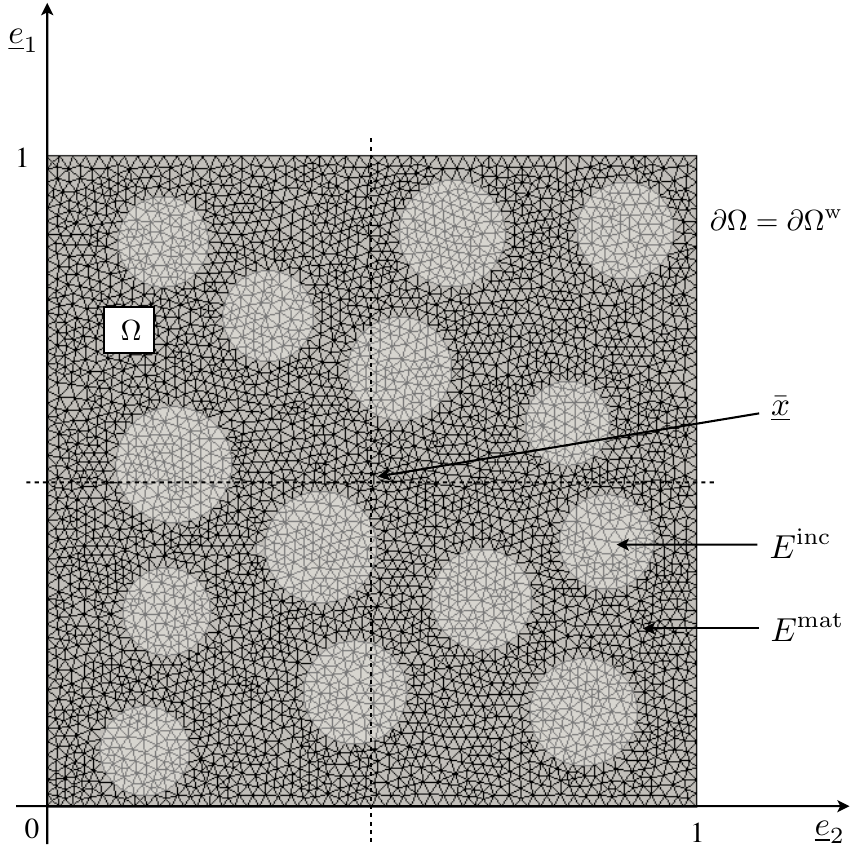}
 \caption{Finite element discretisation of the parametrised SVE problem.}
 \label{fig:Mesh}
\end{figure}

The SVE boundary value problem  \{\eqref{eqref:VirtualWorkparametric},
\eqref{eq:ConstiPara}, \eqref{eq:BCSVE}\}  is parametrised by the material parameters $\Parameter^\text{m}$ and the load parameters $\Parameter^\text{l}$. Domain $\Domain$ is the unit square in 2D, defined by  $\Domain = [0 ,1] \times [0 , 1]$, over which we discretise the elasticity problem with piecewise constant elasticity constants in a conforming way, using linear triangle elements (see figure \ref{fig:Mesh}). Consistently with the previous subsection, the ``truth" finite element SVE problem reads:
\begin{equation} 
\label{eq:FEMproblemExemple}
\begin{array}{l}
\displaystyle  \text{Find} \ \DispFEMHomo(\Parameter) \in \SpaceDispFEMZero \ \text{such that} \ \forall \, \DispTest \in \SpaceDispFEMZero,
\\ \displaystyle  
\int_{\Domain} \StrainBasic(\DispFEMHomo(\Parameter)) : \Hooke(\Parameter) : \StrainBasic(\DispTest) \, \DVolume = - \int_{\Domain} \StrainBasic(\DispFEMLift(\Parameter)) : \Hooke(\Parameter) : \StrainBasic(\DispTest) \, \DVolume \, ,
\end{array} 
\end{equation}
where $\DispFEMLift(\Parameter) = \StrainBasic^\text{M}(\Parameter^\text{l}) (\Pos - \bar{\Pos}) \equiv \bar{\Disp}(\Parameter)$ is the coarse field and $\DispFEMHomo(\Parameter) \equiv \tilde{\Disp}(\Parameter)$ is the fluctuation field in the context of computational homogenisation. It is therefore clear that the linear form associated with the quantity of interest $Q$ is minus the right-hand side of this finite element problem, weighted by the measure of $\Domain$, which is the particular case announced at the end of section \ref{sec:QoI}.

We assume in this example that the material heterogeneities are only parametrised by the elastic contrast $\Parameter^\text{m} \equiv \mu_1 = \frac{\Young^\text{inc}}{\Young^\text{mat}}$, where $\Young^\text{mat}$ is the Young's modulus of the matrix and $\Young^\text{inc}$ is the Young's modulus of the inclusions. The elastic contrast ranges from 0.1 (soft inclusions) to 10 (hard inclusions). The Poisson's ratios of both phases is set to $\Poisson = 0.3$. In this context, the affine representation of the Hooke's elasticity tensor over the parameter domain reads:
\begin{equation}
\forall \, \Parameter \in \ParameterSpace, \,  \forall \, \Pos \in \Domain ,  \quad 
 \Hooke(\Pos,\Parameter) = \Hooke^\text{mat} + (\mu_1-1) \, H^\text{inc}(\Pos) \,  \Hooke^\text{mat} \, .
\end{equation}
In the above equations, function $H^\text{inc}$  is the indicator function of the inclusion phase. It is equal to 1 for a point located in an inclusion and 0 elsewhere. The elasticity tensor of the matrix phase $\Hooke^\text{mat}$ is defined by equation \eqref{eq:Consti}, with $\Young^\text{mat} = 1$, and $\nu^\text{mat} = 0.3$. The affine representation of the compliance tensor over the parameter domain becomes
\begin{equation}
\forall \, \Parameter \in \ParameterSpace, \,  \forall \, \Pos \in \Domain ,  \quad 
 \Compliance(\Pos,\Parameter) = 
 \Compliance^\text{mat} + H^\text{inc}(\Pos) \left(\frac{1}{\mu_1}-1 \right) \Compliance^\text{mat} \, , 
\end{equation}
where $\Compliance^\text{mat}$ is the compliance tensor of the matrix phase.

The independent components of the effective strain constitute parameters $\Parameter^\text{l}$. More precisely, we define $\mu^\text{l}_1 \equiv \mu_2 = \StrainBasic_{11}^\text{M}$, $\mu^\text{l}_2 \equiv \mu_3 = \StrainBasic_{22}^\text{M}$ and $\mu^\text{l}_3 \equiv \mu_4 = \StrainBasic_{12}^\text{M}$. The affine representation of the Dirichlet boundary conditions is 
\begin{equation}
\label{eq:ExpDirichlet}
\forall \, \Parameter \in \ParameterSpace, \,  \forall \, \Pos \in \DDomain ,  \quad 
\DirichletBC(\Pos,\Parameter) = \left( \begin{pmatrix} 1 & 0 \\ 0 & 0 \end{pmatrix} \mu_2 + 
\begin{pmatrix} 0 & 0 \\ 0 & 1 \end{pmatrix} \mu_3 + 
\begin{pmatrix} 0 & 1 \\ 1 & 0 \end{pmatrix} \mu_4 
\right) (\Pos - \bar{\Pos})
\, .
\end{equation}

The parameter domain is restricted to the hypercube $\ParameterSpace = [0.1 , 10] \times [-1,1] \times [-1,1] \times [-1,1]$.



\subsection{Reduced order modelling and certification}

\begin{figure}[htb]
 \centering
 \includegraphics[width=0.6\linewidth]{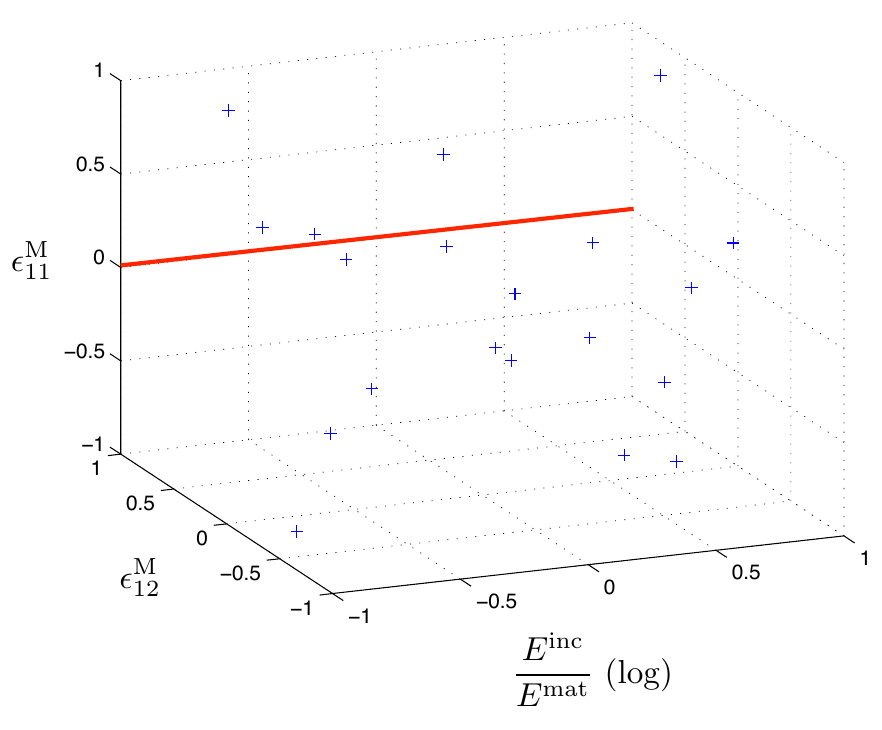}
 \caption{Quasi Monte-Carlo Sampling of the parameter domain. Only three of the four dimensions are represented here. The blue points are the parameter values used to compute the snapshot. The red line is the 1D subdomain over which the ``online" reduced order modelling interpolation and associated error estimation is validated.}
 \label{fig:QMC}
\end{figure}

\paragraph{``Offline" procedure}

In order to deploy the Galerkin-POD methodology, we first sample the parameter domain $\ParameterSpace$. This is done by generating a 4-dimensional Quasi-Random sequence (the Sobol sequence implemented in MATLAB\textsuperscript{\textregistered}) and retaining the first 20 points, which defines $\ParameterSpaceDiscr$ (see figure \ref{fig:QMC}). The corresponding ``truth" finite element displacement and stress fields are stored, for interpolation and error estimation purposes respectively. In our implementation, we store the nodal values of displacements, while we store the stresses at the quadrature points associated with $\SpaceDispFEM$ (\textit{i.e.:} at the centroids of the triangle elements).

Next, we compute the set of lifting functions $(\FunctionAdKinematicAffineI{i})_{i \in \llbracket 1, 3 \rrbracket} \in (\SpaceDispFEM)^3$. In order to do so, we first set $\ParameterZero = (1 \ 0 \ 0 \ 0)^T$ (we therefore make use of the bilinear form associated with an homogeneous material), and then solve the series of corresponding finite element problems (recall equation \eqref{eq:DispLift}) with Dirichlet boundary conditions corresponding to each of the term of affine expansion \eqref{eq:ExpDirichlet}. We can then subtract the lifting $\DispFEMLift(\Parameter)$ from  displacements $\DispFEM(\Parameter)$ for all parameters $\Parameter$ belonging to the sampled parameter domain $\ParameterSpaceDiscr$. A singular value decomposition is performed on the resulting bubble fields (they vanish over $\DDomain \equiv \DDomainDirichlet$). The first $\NumberReducedBasis$ spatial modes of this proper orthogonal decomposition are then stored (we store their nodal values).

The lifting for the stress needs not be performed as $\VolumeForce(\Pos) =\V{0} ,\ \forall \Pos \in \Domain$ and $\DDomainNeumann = \{ \}$. A Snapshot POD in the sense of the energy norm associated to $\ParameterZero$ is performed on the set of sampled finite element stress, and the first $\NumberReducedBasisStress$ spatial modes are stored (we store their values at the quadrature points).

Last, we can pre-compute all the operators that allow for the ``offline"/``online" split of the numerical complexity associated to interpolation and verification as described througout sections \ref{sec:MOR} and \ref{sec:CRE}.

\paragraph{``Online" procedure}

In the online phase, for a particular $\Parameter \in \ParameterSpace$ of interest, we need to evaluate the lifting functions $\DispFEMLift(\Parameter)$ and $\StressFEMLift(\Parameter)$ (the latter is null in our case), and interpolate the complementary parts of the displacement and stress fields, $\DispReducedHomo(\Parameter)$ and $\StressFEMReducedHomo(\Parameter)$ respectively. The first set of these operations is done by direct evaluation (equation \eqref{eq:DispLift}), while the latter requires the solution of the projected problems for the displacement (equation \eqref{eq:GalerkinReducedBasis}), for the hierarchically enriched displacement (equation \eqref{eq:GalerkinReducedBasisEnriched}) and for the equilibrated stress in the finite element sense (equation \eqref{eq:MiniStress}). 

Finally, the upper bound for the ``online" error associated with the Galerkin-POD is obtained by comparison of the recovered stress field $\StressHat(\Parameter)$ to the stress field $\StressReduced(\Parameter)$ obtained by applying the constitutive relation to the reduced approximation of the displacement $\DispReduced(\Parameter)$, as explained in section \ref{sec:CRE}. The lower bound is obtained by comparing the hierarchically enriched displacement field $\DispReducedEnhanced(\Parameter)$ to the reduced displacement $\DispReduced(\Parameter)$.

\paragraph{Remark on the relationship between error of homogenisation and error in energy norm: } we come back to the statements made in section \ref{sec:QoI} about the particular form of quantity of interest that is used in the context of our computational homogenisation scheme computational for isotropic composite materials. If the effective strain that is applied online is the pure shear strain 
$ \StrainBasic^\text{M}(\Parameter^\text{l}) = \frac{1}{2} \left( \UnitVector{1} \otimes \UnitVector{2} +  \UnitVector{2} \otimes \UnitVector{1} \right)$, then the quantity of interest $Q(\Parameter)$ extracts the effective shear modulus $G^\text{M}(\Parameter^\text{m})$. The error in this constant that introduced by reduced order modelling is directly minus the square of the energy norm of the error field, weighted by the measure of the domain. Therefore, we have 
\begin{equation}
 - \frac{(\nu^\text{low}_\text{G})^2}{| \Domain|} \geq G^\text{M,h} - G^\text{M,r} \geq - \frac{ (\nu^\text{up}_\text{G})^2}{| \Domain|} \, ,
\end{equation}
where $G^\text{M,h}$ is the effective shear modulus obtained using the finite element method directly, while $G^\text{M,r}$ is its reduced order modelling approximation, $\nu^\text{low}_\text{G} \equiv \LowerBound(\Parameter)$ and $\nu^\text{up}_\text{G} \equiv \CRE(\Parameter)$. The same concept can be applied with $ \StrainBasic^\text{M}(\Parameter^\text{l}) = \UnitVector{1} \otimes \UnitVector{1} $ to obtain bounds on the effective first Lam\'e constant that make direct use of the energy norm of the reduced order modelling error in both numerical tests. After algebraic manipulations, this reads 
 \begin{equation}
- \frac{(\nu^\text{low}_{\lambda + \text{2G}})^2}{| \Domain|} 
+ 2 \, \frac{(\nu^\text{up}_\text{G})^2}{| \Domain|} 
\geq \lambda^\text{M,h} - \lambda^\text{M,r} \geq - \frac{ (\nu^\text{up}_{\lambda+ \text{2G}})^2}{| \Domain|} 
+ 2 \, \frac{(\nu^\text{low}_\text{G})^2}{| \Domain|}  \, ,
\end{equation}
where $\nu^\text{up}_{\lambda + \text{2G}} \equiv \CRE(\Parameter)$ in this second setting,  and $\nu^\text{low}_{\lambda + \text{2G}} \equiv \LowerBound(\Parameter)$. It is important to notice that these bounds only take into account the reduced order modelling error. The error due to the finite element discretisation, and more importantly the error due to using a material sample of finite size in this simplified computational homogenisation framework, are both ignored.


\subsection{Numerical results}

We validate our error bounding methodology on a subdomain of the parameter domain $\mathcal{P}^\text{eval} = \{ \Parameter \in \ParameterSpace \, | \,  \mu_2=0, \,  \mu_3=0 \ \text{and} \ \mu_4 = 1 \}$, which is illustrated in figure \ref{fig:QMC} and corresponds to applying a shear load to the SVE. Notice that none of the sampled parameter values actually belong to this subdomain. The numerical results are given in figures \ref{fig:Effectivity} and \ref{fig:Convergence} and are commented below. In order to ease the interpretation of the numerical results, we will use the relative error bounds $\CRE^\text{,rel}$ and $\LowerBound^\text{,rel}$ defined by
\begin{equation} 
\label{eq:RelativeCRE}
\forall \, \Parameter \in \ParameterSpace , \quad \CRE^\text{,rel}(\Parameter) \Def
\frac{ \displaystyle 
\CRE(\Parameter)  
}{  \displaystyle 
\| \DispReducedEnhanced(\Parameter) \|_{\Hooke(\Parameter)}  
}
\quad \text{and} \quad
\LowerBound^\text{,rel}(\Parameter) =
\frac{ \displaystyle
\LowerBound(\Parameter)
}{ \displaystyle
\| \DispReducedEnhanced(\Parameter) \|_{\Hooke(\Parameter)}  
} \, ,
\end{equation}
which satisfy the bounding properties
\begin{equation}
\forall \, \Parameter \in \ParameterSpace , \quad  
\LowerBound^\text{,rel}(\Parameter)
\leq
\frac{ \displaystyle \| \ErrorReduced(\Parameter)  \|_{\Hooke(\Parameter)}}{  \displaystyle \| \DispReducedEnhanced(\Parameter) \|_{\Hooke(\Parameter)}   } \leq \CRE^\text{,rel}(\Parameter) \, .
\end{equation}

\subsubsection{Effectivity of the error bounds}

\begin{figure}[htb]
 \centering
 \includegraphics[width=1\linewidth]{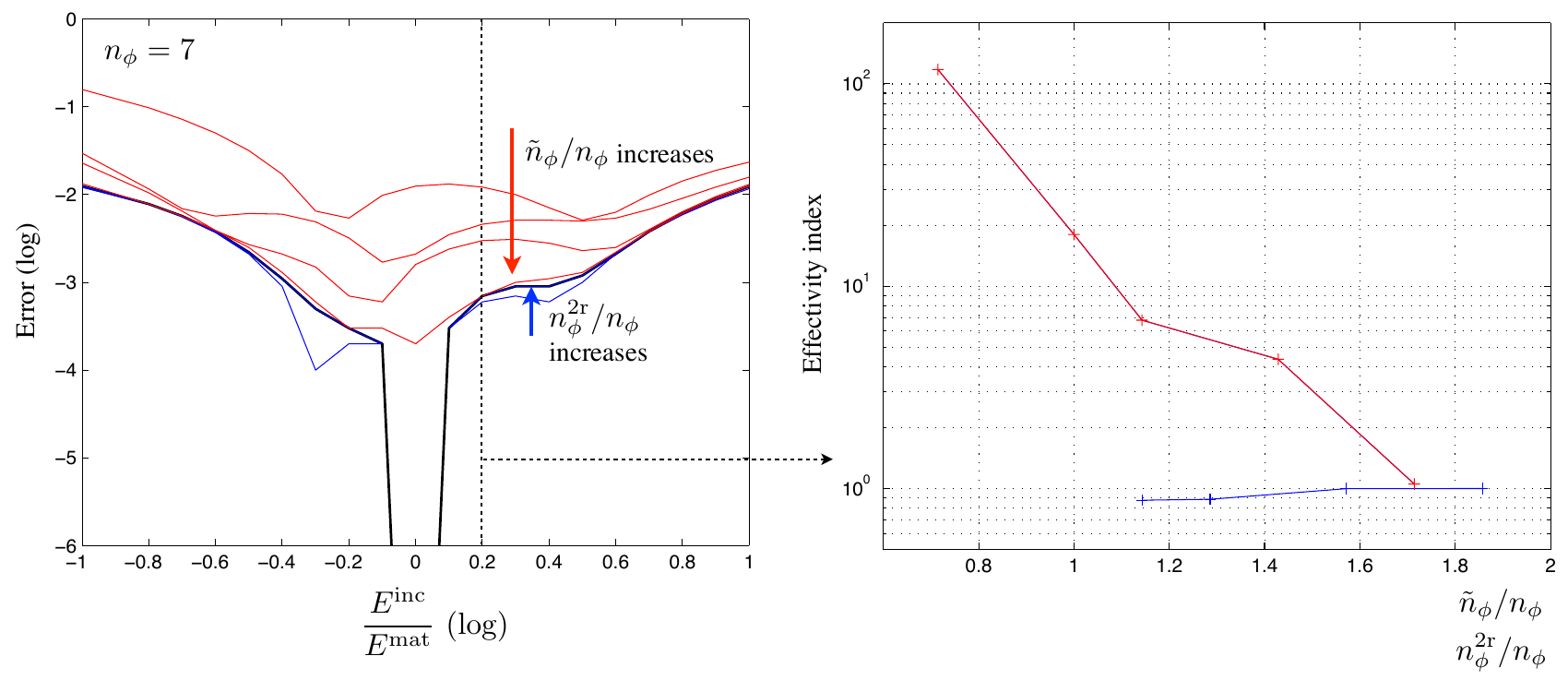}
 \caption{Effectivity of the upper and lower bounds for the ``online" error associated to the reduced order modelling approach as a function of the elastic contrast of the composite structure and of the expansion order of the two surrogates used to obtain the bounds.}
 \label{fig:Effectivity}
\end{figure}

First, we show in figure \ref{fig:Effectivity} the influence of parameters $\NumberReducedBasisStress$ and $\NumberReducedBasisEnriched$ on the effectivity of the upper and lower bounds for the error introduced by the reduced order modelling approximation. 

For demonstration purposes, $\NumberReducedBasis$ is arbitrarily set to 7. The left-hand side graph in figure \ref{fig:Effectivity} shows the evolution of the exact relative error in energy norm $ \frac{ \| \ErrorReduced(\Parameter) \|_{\Hooke(\Parameter)} }{   \| \DispReducedEnhanced(\Parameter) \|_{\Hooke(\Parameter) } }$, the upper bound $\nu^\text{up,rel}$ for this error measure and the lower bound $\nu^\text{low,rel}$ as a function of the elastic contrast $\mu_1$. The upper bound is first obtained by using an expansion of order 5 for the stress surrogate (\textit{i.e.} $\NumberReducedBasisStress=5$). This number is then increased to 12 in an incremental manner. Similarly, the lower bound is first obtained using an expansion of order 8 for the enriched displacement surrogate (\textit{i.e.} $\NumberReducedBasisEnriched=8$), which is then increased to 13 in an incremental manner. It clearly appears that both the lower and upper bounds consistently sharpen with increasing dimensions of the corresponding reduced spaces.

The sharp decrease of the exact error around point $(1 \ 0 \ 0 \ 1 )^T$ of the parameter domain can be explained by the fact that our Galerkin-POD model provides the exact finite element solution at this particular point. Indeed, the exact solution is obtained by direct evaluation of the lifting $\DispFEMLift(\Parameter)$, due to the choice that we made for $\ParameterZero$. The complementary part of the reduced solution $\DispReducedHomo(\Parameter)$ vanishes at that point.

The effectivity results stated previously are emphasised by the right-hand side graph in figure \ref{fig:Effectivity}. The elastic contrast is fixed at a value of $1.6$. The effectivity index of the two bounds, defined by
\begin{equation}
\forall \, \Parameter \in \ParameterSpace, \quad 
\left\{ \begin{array}{l}
\displaystyle
\theta^\text{up} = \frac{\nu^\text{up}}{
 \| \ErrorReduced(\Parameter) \|_{\Hooke(\Parameter)} 
} 
\\ \displaystyle
\theta^\text{low} = \frac{\nu^\text{low}} 
{  
\| \ErrorReduced(\Parameter) \|_{\Hooke(\Parameter)}  
}   \, ,
\end{array} \right.
\end{equation}
are plotted as a function of the ratio between the the dimension of the corresponding reduced spaces, $\NumberReducedBasisStress$ and $\NumberReducedBasisEnriched$, and the dimension of the reduced space used to compute the initial approximation of the solution (\textit{i.e.} $\NumberReducedBasis=7$). Both bounds seem to converge quickly to an effectivity of 1. The lower bound is extremely effective even when using a very small number of additional POD modes. The upper bound needs slightly more computational effort to reach the expected numerical efficiency. However, the sharpness of the upper bound could surely be controlled in an adaptive manner using the available distance between the upper and the lower bounds.

\paragraph{Remark:} It should be noticed that in spite of the fact that the bounds seem to converge to the exact error with increasing dimensions of the corresponding reduced spaces, we have no reason to believe that this behaviour is asymptotically true. Indeed, the information that is available to train the surrogate models is obtained by sampling the parameter domain. Therefore, the exact finite element stress and exact finite element displacement cannot, in general, be obtained by simply increasing $\NumberReducedBasisStress$ and $\NumberReducedBasisEnriched$. In fact, we should observe a stagnation of the effectivities when using reduced spaces of large dimensions to compute the bounds. To alleviate this potential issue, the only solution is to refine the sampling of the parameter domain in an adaptive manner.

\subsubsection{Convergence of the error bounds}

\begin{figure}[htb]
 \centering
 \includegraphics[width=1\linewidth]{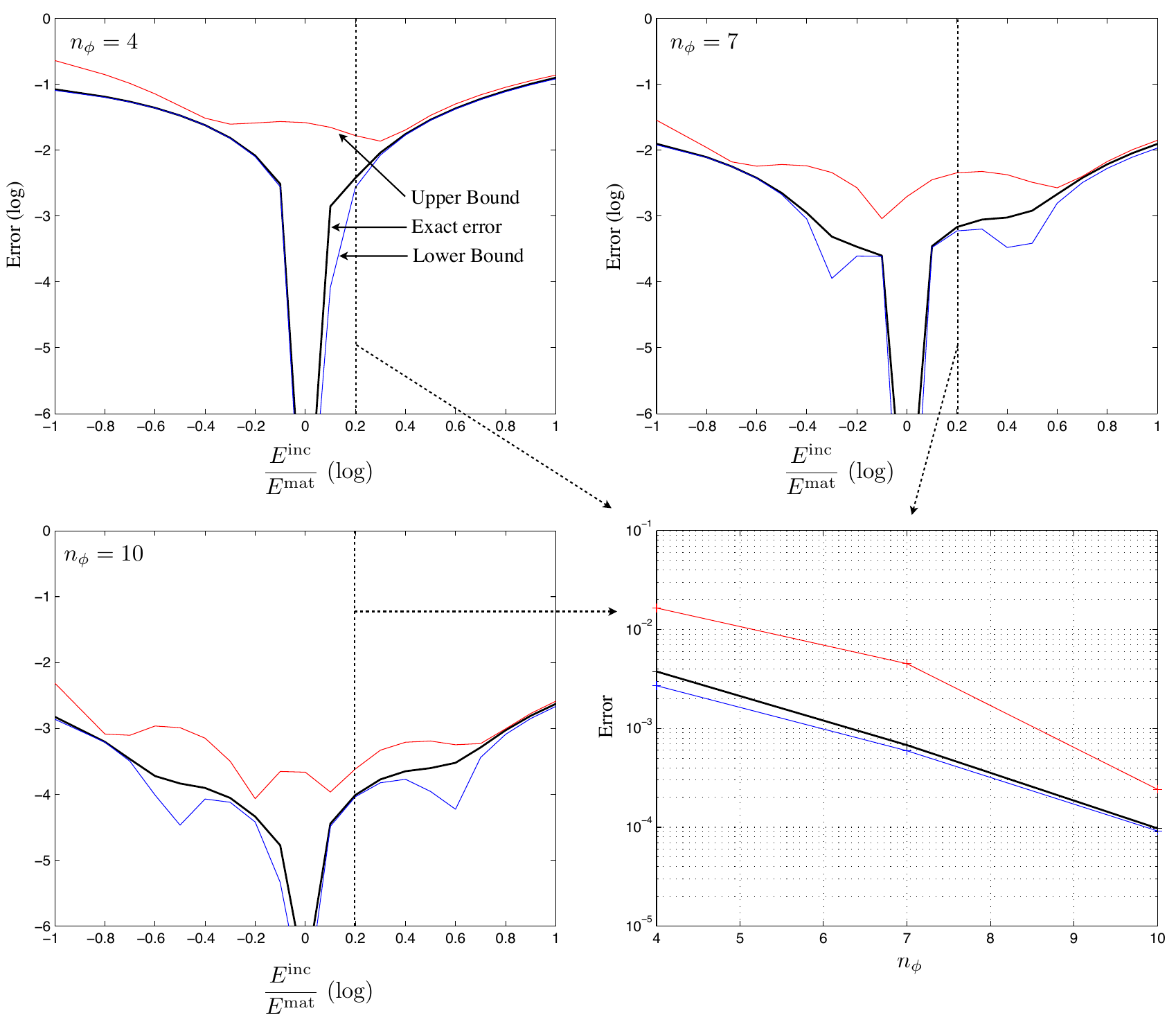}
 \caption{Convergence of the upper and lower bounds of the reduced order modelling error as a function of the order of the POD expansion used to approximate the solution of the parametrised problem of computational homogenisation over the parameter domain.}
 \label{fig:Convergence}
\end{figure}

Next, we show how the proposed bounds converge with the error in the reduced order modelling approximation. To do so, we set $\NumberReducedBasisStress = \NumberReducedBasis+2$ and $\NumberReducedBasisEnriched = \NumberReducedBasis+1$. In this manner, the numerical parameters of the lower and upper bounds are constrained to follow the order of the surrogate model that we want to certify. In figure \ref{fig:Convergence}, we show the convergence of the exact relative error $ \frac{ \| \ErrorReduced(\Parameter) \|_{\Hooke(\Parameter)} }{   \| \DispReducedEnhanced(\Parameter) \|_{\Hooke(\Parameter) } }$ and that of the upper and lower bounds $\nu^\text{up,rel}$ and $\nu^\text{low,rel}$ as a function of $\NumberReducedBasis$. As for the previous results, the plots are given for the subdomain $\mathcal{P}^\text{eval}$. 

In the right-hand side bottom graph, we fix $\mu_1=1.6$ to ease the interpretation. These partial results show that the bounds have approximately the same convergence rate as the reduced order modelling approximation that we want to certify.

\section{Discussion and conclusion}

We have presented a novel and conceptually simple way to build lower and upper bounds for the error arising in POD-based reduced order models for parametrised elasticity problems. The upper bound relies on the error in the constitutive relation, which requires to construct and evaluate a reduced order model for the stress field. This bound does not rely on the Galerkin orthogonality, which potentially allows us to apply it to other type of reduced order modelling techniques (Kriging, Moving Least-Squares approximations, response surface method, \textit{etc.}). We have shown good numerical convergence and effectivity properties of both bounds, which are the basic requirements to perform adaptivity.

In terms of implementation, we have shown that for each step of the bounding procedure, and adequate ``offline"/``online" split of the numerical complexity guarantees that the certification remains efficient numerically. However, the intrinsic nature of projection-based reduced order modelling makes the implementation complex, and highly intrusive.

The bounding techniques have been applied to a basic example of multiscale modelling, for which the requirements in terms of cost reductions are still rather stringent. On the way, we have laid some foundations for future expansion of this work to more complex microstructures and micro/macro relationships. \\

Our next step is to directly evaluate the surrogate error in terms of quantities of interest. In the case of elasticity, this is a rather straightforward step, as the goal-oriented error estimation techniques based on duality are very well-established. From there, we can imagine developing an adaptive reduced order modelling approach. To do that, efficient procedures need to be devised to control the dimension of the reduced space used to approximate the solution, the sharpness of the error bounds and the sampling of the parameter domain.

In terms of homogenisation, the goal-oriented certification will require to devise two type of error bounds, one for the averaged information, controlled by a global quantity of interest, and one for the local information at the material level. The local quantities of interest are very often the maximum value of some stress or strain components. Unfortunately, the corresponding dual problems are not reducible in the general case. An approximation will be necessary, at the probable cost of loosing the bounding property, or not bounding the right quantity.

At last, we hope that the proposed approach complements existing work in the area of error estimation for reduced order modelling. A formal and numerical comparison of the methods available will be necessary. More importantly, we believe that this contribution will help progressing towards the certification of non-affine and nonlinear parametrised boundary value problem, for instance by using the classical extensions of the the constitutive relation error to nonlinear problems.


\section*{Acknowledgments}

The authors acknowledge the financial support of EPSRC under grant EP/J01947X/1: \textit{Towards rationalised computational expense for simulating fracture over multiple scales}(RationalMSFrac) and the support of the European Research Council under Starting Independent Research Grant agreement No. 279578: \textit{Towards real time multiscale simulation of cutting in non-linear materials with applications to surgical simulation and computer guided surgery} (RealTCut).
 

\bibliographystyle{plain}
\bibliography{bibliography.bib}

\end{document}